\documentclass[a4paper,11pt,twoside,fleqn,article]{memoir}
\usepackage[T1]{fontenc}
\usepackage[latin9]{inputenc}
\usepackage{lmodern}

\newcommand{\ifaq}[1]{#1}
\newcommand{\ifac}[1]{}

\fixpdflayout

\usepackage{graphics}
\usepackage[usenames,dvipsnames]{xcolor}

\input{bar.tex}
\usepackage{pifont}
\usepackage[final]{microtype}
\usepackage[british,french]{babel}
\usepackage{graphicx}
\usepackage{cancel}
\DeclareMathSizes{9}{8}{5}{5}
\DeclareMathSizes{10}{10}{7}{5}
\DeclareMathSizes{10.95}{10}{7}{5}

\addto\extrasfrench{\providecommand{\fg}{\ifdim\lastskip>\z@\unskip\fi~\frqq}}
\usepackage{bm}

\usepackage[geometry]{ifsym}

\newlength{\currentparskip}

\definecolor{gristc}{rgb}{0.85,0.85,0.85}
\newcounter{TODOS}

\newlength{\largeurpage}
\newcommand{\NumT}{\smallskip\par\noindent\colorbox{gristc}{\begin{minipage}{\largeurpage}\hfill{\normalsize\sffamily\bfseries\selectfont{}~\#Todo = \arabic{TODOS}~}\mbox{}\end{minipage}}\newline\smallskip

}
\renewcommand{\NumT}{}

\captionstyle[\raggedright]{}
\newcommand{\captionsymb}{\captionsymba}
\newcommand{\captionsymba}{\raisebox{-2.0pt}{\mbox{\color{grisc}\!\FilledSmallCircle\:}}%
\renewcommand{\captionsymb}{\captionsymbb}%
}
\newcommand{\captionsymbb}{\renewcommand{\captionsymb}{\captionsymba}}

\captionnamefont{\captionsymb\small\sffamily\bfseries\boldmath\color{grisf}}
\captiontitlefont{\small}

\makeatletter
\newlength{\figrulesep}
\newcommand{\topfigrule}{\kern-0.8pt%
\noindent\raggedright{\color{grisf}\rule[-\figrulesep]{2in}{0.8pt}}}

\makeatother

\usepackage[thref,thmmarks]{ntheorem}
\theoremstyle{plain}
\theoremheaderfont{\hspace*{1pt}\raisebox{-2.0pt}{\mbox{\color{grisc}\FilledSmallTriangleRight}}~\upshape\sffamily\bfseries\boldmath\color{grisf}}
\theoremstyle{break}
\theorembodyfont{\normalfont}
\theoremseparator{\vspace{3pt}}
\theoremsymbol{}
\newcounter{lem}
\newcounter{thm}
\newcounter{def}

\newtheorem{thmen}[thm]{Theorem}
\newtheorem{charen}[thm]{Characterization}

\newtheorem{coren}[thm]{Corollary}
\theoremsymbol{\hfill\raisebox{-2.0pt}{\mbox{\color{grisc}\FilledSmallTriangleLeft}}}

\newtheorem{charennp}[thm]{Characterization}
\theoremsymbol{}

\theoremstyle{plain}
\theoremsymbol{}
\theorembodyfont{\normalfont}
\theoremheaderfont{\hspace*{1pt}\raisebox{-2.0pt}{\mbox{\color{grisc}\FilledSmallTriangleRight}}~\sffamily\itshape\bfseries\boldmath\color{grisf}}

\newtheorem{lemen}[lem]{Lemma}
\theoremsymbol{\hfill\raisebox{-2.0pt}{\mbox{\color{grisc}\FilledSmallTriangleLeft}}}

\newtheorem{lemennp}[lem]{Lemma}
\theoremsymbol{}

\theorembodyfont{\normalfont}
\theoremheaderfont{\vspace*{1pt}\sffamily\bfseries\color{grisf}\boldmath}

\newtheorem{defen}[def]{Definition}
\theoremsymbol{}

\theoremheaderfont{\sffamily\bfseries\boldmath\slshape\selectfont\color{grisf}}
\newcounter{exrem}
\newtheorem{remen}[exrem]{Remark}

\newtheorem{exen}[exrem]{Example}

\theoremsymbol{\hfill\raisebox{-2.0pt}{\mbox{\color{grisc}\FilledSmallTriangleLeft}}}
\theoremstyle{nonumberplain}

\theorembodyfont{\normalfont}
\theoremheaderfont{\sffamily\bfseries\boldmath\selectfont\color{grisc}}
\theoremseparator{.~}
\newcounter{proof}
\newtheorem{proofen}[proof]{Proof}

\theoremsymbol{}
\newtheorem{proofennoend}[proof]{Proof}

\usepackage{array}
\usepackage{multirow} 
\usepackage{multicol}
\usepackage{verbatim}
\usepackage{color}
\usepackage[color,matrix,arrow,all]{xy}

\newcommand*{\mysection}[2][\mysectionname]{%
	\def\mysectionname{#2}%
	\subsection*{#2}%
}

\usepackage{quoting}
\quotingsetup{vskip=4pt}

\usepackage{listings}
\input{headers.tex}

\let\oldwedge\wedge\renewcommand{\wedge}{\!\oldwedge\!}
\let\oldvee\vee\renewcommand{\vee}{\!\oldvee\!}
\let\oldbullet\bullet\renewcommand{\bullet}{{\color{grisf}\oldbullet}}\newcommand{\puce}{\oldbullet}

\newcommand{\N}{\mathcal{N}}
\newcommand{\M}{\mathcal{M}}
\ifaq{}
\ifac{}
\ifaq{}
\ifac{}

\newcommand{\Hy}{\mathcal{H}}

\newcommand{\Mnp}{\mathrm{M}}
\newcommand{\Mp}[1]{\Mnp\left(#1\right)}

\newcommand{\Vertexp}[1]{\mathrm{V}\!\left(#1\right)}

\usepackage[tight]{shorttoc}

\makeatletter
\ifac{%
\usepackage{ifthen}
\newboolean{@multitoc@toc}
\newboolean{@multitoc@lof}
\newboolean{@multitoc@lot}
\newboolean{@multitoc@loa}
\setboolean{@multitoc@toc}{true}
\setboolean{@multitoc@lof}{true}
\setboolean{@multitoc@lot}{true}
\setboolean{@multitoc@loa}{true}
\newcommand*{\multicolumntoc}{2}%
\newcommand*{\multicolumnlof}{2}%
\newcommand*{\multicolumnlot}{2}%
\newcommand*{\multicolumnloa}{2}%
\let\@multitoc@starttoc\@starttoc%
\renewcommand*{\@starttoc}[1]{%
	\ifthenelse{\boolean{@multitoc@toc}\and\equal{#1}{toc}}{%
	\begin{multicols}{\multicolumntoc}%
	\@multitoc@starttoc{#1}%
	\end{multicols}%
	}{}%
	\ifthenelse{\boolean{@multitoc@lot}\and\equal{#1}{lot}}{%
	\begin{multicols}{\multicolumnlot}%
	\@multitoc@starttoc{#1}%
	\end{multicols}%
	}{}%
	\ifthenelse{\boolean{@multitoc@lof}\and\equal{#1}{lof}}{%
	\begin{multicols}{\multicolumnlof}%
	\@multitoc@starttoc{#1}%
	\end{multicols}%
	}{}%
	\ifthenelse{\boolean{@multitoc@loa}\and\equal{#1}{loa}}{%
	\begin{multicols}{\multicolumnloa}%
	\@multitoc@starttoc{#1}%
	\end{multicols}%
	}{}%
}%
}
\makeatother

\maxsecnumdepth{subsubsection}
\maxtocdepth{subsubsection}

\usepackage[tight]{minitoc}
\nomtcrule

\renewcommand{\footnoterule}{%
	{\color{grisf}%
   \kern -2.8pt   \hrule width 2in height 0.8pt   \kern 2pt }%
}
\renewcommand*{\hrulefill}[1][1pt]{\leavevmode\leaders\hrule height #1 \hfill \kern 0pt}

\definecolor{grisc}{rgb}{0.6,0.6,0.6}
\definecolor{grisf}{rgb}{0.4,0.4,0.4}
\definecolor{noir}{RGB}{0,0,0}
\definecolor{black}{RGB}{0,0,0}
\usepackage{colortbl}
\arrayrulecolor{grisc}
\usepackage[vlined,french,norelsize]{algorithm2e}

     {\begin{list}{}%
             {\setlength{\leftmargin}{#1}}%
             \item[]%
     }
     {\end{list}}

\SetAlCapSty{mysty}
\DontPrintSemicolon

\SetAlgoCaptionLayout{myalgocaplay}
\SetAlCapFnt{\sffamily\bfseries\selectfont\color{grisf}}

\SetCommentSty{Stylecomment}{}{}
\SetNlSty{Stylenombres}{}{}
\SetNlSkip{0.4em}
\SetInd{0.5em}{1.5em}
\SetAlFnt{\small}
\SetAlCapNameFnt{\small}
\SetArgSty{}
\SetKwSty{paperxxvKwSty}

\newcommand{\AlgoEnglish}{
	\renewcommand{\algorithmcfname}{Algorithm}%
	\renewcommand{\algorithmautorefname}{algorithm}%
	\SetAlgoCaptionSeparator{:}
	\SetKw{Yield}{yield}
	\SetKwIF{If}{ElseIf}{Else}{if}{{\!}:}{elsif}{else}{endif}
	\SetKwFor{For}{for}{{\!}:}{endfor}
	\SetKwFor{Foreach}{foreach}{{\!}:}{endfor}
	\SetKwFor{While}{while}{{\!}:}{endw}
	\SetKwBlock{Begin}{:}{end}
	\SetKw{Return}{return}
	\SetKw{ReturnVoid}{return}
	\SetKw{Write}{write}
}

\firmlists*

\makeatletter
\renewcommand{\algocf@Vline}[1]{%
                \strut\par\nointerlineskip
                \algocf@push{\skiprule}%
                \hbox{%
                        {\color{grisc}\vrule width 0.4pt}
                        \vtop{%
                                \algocf@push{\skiptext}%
                                \vtop{\algocf@addskiptotal\advance\hsize by -\skiplength #1}%
                                {\color{grisc}\Hlne}
                        }%
                }\vskip\skiphlne%
                \algocf@pop{\skiprule}%
                \nointerlineskip%
        }

        \renewcommand{\algocf@Vsline}[1]{%
                \strut\par\nointerlineskip%
                \algocf@push{\skiprule}%
                \hbox{%
                        {\color{grisc}\vrule width 0.4pt}
                        \vtop{\algocf@push{\skiptext}%
                        \vtop{\algocf@addskiptotal\advance\hsize by -\skiplength #1}}%
                }%
                \algocf@pop{\skiprule}%
        }
\makeatother

\let\stdtoc\tableofcontents
\renewcommand*\tableofcontents{{%
\renewcommand*\MakeUppercase[1]{##1}\stdtoc}}
\nouppercaseheads

\makeatletter
\renewcommand{\dotfill}{\leavevmode\xleaders\hbox{\hspace*{0.25em}{\small.}\hspace*{0.25em}}\hfill\kern0pt}
\makeatother

\usepackage{enumitem}
\setlist{itemsep=0ex,partopsep=0ex,topsep=0ex,parsep=0ex,leftmargin=\parindent,labelindent=\parindent,font=\sffamily\bfseries\selectfont\color{grisc}\boldmath}
\setdescription{labelindent=0.5\parindent}
\setenumerate{leftmargin=\parindent,labelindent=\parindent,font=\sffamily\bfseries\boldmath\color{grisc},label=\arabic*}

\newcommand{\englishlabel}{~\color{grisc}$\puce$}

\newcommand{\englishlist}{\renewcommand{\labelitemi}{\englishlabel}\renewcommand{\labelitemii}{\englishlabel}}

\newcommand{\English}{
	\captiondelim{:{ }}
	\AlgoEnglish%
	\selectlanguage{british}%
	\englishlist
}

\def\clap#1{\hbox to 0pt{\hss#1\hss}}

\usepackage{amsfonts}
\usepackage{amssymb}\usepackage{amsmath}
\usepackage[fleqn]{nccmath}

\renewcommand{\ldots}{\hspace{-1pt}...}
\ifac{\usepackage{sfmath}}
\newcommand{\mysymbolr}{\color{grisc}\right.\color[RGB]{0,0,0}}
\newcommand{\mysymbol}{\color{grisc}\left|}
\newcommand{\mydisplaystyle}[1]{\vspace{-3pt}\begin{align*}~\mysymbol~\color[RGB]{0,0,0}\hfill\vcenter{\vspace{1.5pt}#1}\hfill\mysymbolr\end{align*}}

\let\oldsmallsetminus\smallsetminus\renewcommand{\smallsetminus}{\!\oldsmallsetminus\!}
\begin{document}
\ifaq{
}
\let\oldabovedisplayskip\abovedisplayskip\setlength{\abovedisplayskip}{0.25\oldabovedisplayskip}
\let\oldbelowdisplayskip\belowdisplayskip\setlength{\belowdisplayskip}{0.25\oldbelowdisplayskip}
\let\oldabovedisplayshortskip\abovedisplayshortskip\setlength{\abovedisplayshortskip}{0.25\oldabovedisplayshortskip}
\let\oldbelowdisplayshortskip\belowdisplayshortskip\setlength{\belowdisplayshortskip}{0.25\oldbelowdisplayshortskip}
%
\renewcommand{\thesection}{\arabic{section}}
\bibliographystyle{alpha}
\bibdata{}

\let\oldleq\leq\renewcommand{\leq}{\!\oldleq\!}
\newcommand{\less}{<}
\let\oldless\less\renewcommand{\less}{\!\oldless\!}

\English
\chaptermarkS{Hypergraph Acyclicity Revisited}%
\thispagestyle{empty}
\vspace*{\beforechapskip}
\chapter*{Hypergraph Acyclicity Revisited}
{\small\sffamily\bfseries 
\selectfont
\noindent \color{grisf}Johann Brault-Baron --- LSV/ENS--Cachan/Inria --- \today}
\vspace*{\afterchapskip}
\paragraph{Abstract}
The notion of graph acyclicity has been extended
to several different notions of hypergraph acyclicity,
 in increasing order of generality: \emph{gamma} acyclicity,
\emph{beta} acyclicity, and \emph{alpha} acyclicity,
that have met a great interest in many fields.

We prove the equivalence between the numerous characterizations
of each notion with a new, simpler proof, in a self-contained manner.
For that purpose, we introduce new notions of alpha, beta and gamma leaf that 
allow to define new ``rule-based'' characterizations of each notion.

The combined presentation of the notions is completed with
a study of their respective closure properties.
New closure results are established, and alpha, beta and gamma
acyclicity are proved optimal w.r.t. their closure properties.

\mysection{Introduction}
\noindent
The notion of graph acyclicity has been extended
to several different ``degrees of acyclicity'' of hypergraphs.
One can cite, in increasing order of generality: \emph{gamma} acyclicity,
\emph{beta} acyclicity, and \emph{alpha} acyclicity.
Each of these notions admits many different characterizations,
and has found applications in database theory \cite{desi}
(see also \cite{PapadimitriouY99,logcfl,hypertree,FlumFG02} 
for examples of application of alpha acyclicity in this context), 
constraint satisfaction problems \cite{rina,ordyniak_et_al:LIPIcs:2010:2855},
and finite model theory \cite{BaganDurandGrandjean07},
and they are ``also of interest from a purely graph-theoretic viewpoint'' \cite{Fagin83degreesof}.
The present author's paper \cite{moi:csl} is a simple, canonical example of application in finite model theory; we develop on this below.

The different characterizations of gamma, beta and alpha acyclicity were presented (even introduced, for beta and gamma) 
jointly in \cite{Fagin83degreesof},
but are mixed with a great number of characterizations that focus on database aspects.
Furthermore, the proofs are not all self-contained, and often rely on non-trivial graph theoretic results.
We argue that these notions are interesting by themselves, independently of database concerns,
and that proving the equivalence of the different characterizations 
can be done in a self-contained fashion, and that it is even easier to do so.
In addition to that, this is an opportunity to incorporate some new characterizations
that have been proved since \cite{Fagin83degreesof}, and to take advantage of 
the simplifying effects of the framework introduced in \cite{theseDuris,betagamma}.

This paper aims at providing the main characterizations of gamma, beta and alpha
acyclicity, in a structured framework, in a (hyper)graph theoretic perspective.
It takes as a starting point that the notion of acyclic graph admits two very different characterizations:
\begin{description}
\item[(1)] The graph does not contain a cycle.
\item[(2)] The graph is a forest, that is to say:
 it can be reduced to the empty graph by repeatedly removing \emph{leaves}, i.e. vertices that have at most one neighbour.
\end{description}
We observe that each of the three mentioned degrees of acyclicity admit two main classes of characterizations,
that roughly correspond to generalizations of these ones.
As an example, beta acyclicity admits the two following characterizations:
\begin{description}
\item[(1)] The hypergraph does not contain a \emph{beta} cycle.
\item[(2)] One can reduce the hypergraph to the empty hypergraph by removing repeatedly \emph{nest points} (or \emph{beta} leaves).
\end{description}
This illustrates that each characterization of some type of acyclicity is either:
\begin{description}
\item[(1)] a global property (the absence of a certain kind of ``cycle''), or
\item[(2)] a property of reducibility to the empty hypergraph through a certain ``reduction'' process.
\end{description}
These two families of characterizations find different practical use.
In finite model theory, for certain classes of queries, we can make a connection between the complexity of a query
and its associated hypergraph.
The characterizations of type (1) are useful to prove hardness results 
(the presence of a cycle makes a query ``hard''), while
those of type (2) are interesting to prove easiness results
(the reduction process provides a way to reduce repeatedly the problem up to a trivial one).
As an illustration, \cite{BaganDurandGrandjean07} and our paper \cite{moi:csl} faithfully follow this pattern, in the case of alpha and beta acyclicity respectively.

An observation that motivates this work is 
that the ``hard part'' of the proof of the equivalence between the different characterizations lies in fact in the proof that some characterization
of type (1) implies some characterization of type (2).
In the case of the given example, the implication (1)$\Rightarrow$(2)
relies on a result by Brouwer and Kolen \cite{nest}
that a beta acyclic (defined as a characterization of type (1)) 
hypergraph has a nest point, that we call \emph{beta leaf}.

Our main contribution is twofold. The technical one is a result
stronger than the result of Brouwer and Kolen, 
but yet having a more natural and shorter proof than the one in \cite{nest}.

The second aspect of the main contribution is the introduction of two notions similar to that of \emph{nest point} (or: \emph{beta leaf}),
that are called respectively \emph{alpha leaf} and \emph{gamma leaf},~ such that the following are equivalent:
\begin{description}
\item[(1)] The hypergraph does not contain a \emph{alpha} (resp. \emph{gamma}) cycle.
\item[(2)] We can reduce the hypergraph to the empty hypergraph by repeatedly removing \emph{alpha} (resp. \emph{gamma}) leaves.
\end{description}
These new characterizations fit the ``rule-based characterizations'' framework of \cite{theseDuris,betagamma},
and have the advantage over known rule-based characterizations that they actually consist of a \emph{single} rule, 
that is easy to make deterministic, that is: remove \emph{all} the alpha (resp. beta, gamma) leaves. 

By using the same natural idea as in the case of beta acyclicity, 
it is once again proved that (1)$\Rightarrow$(2), which is the hard part.
As a side benefit, this proof is fully self-contained,
by contrast with the proof (1)$\Rightarrow$(2) in \cite{desi}, that relies on a graph theoretic result,
proved in \cite{Golumbic} (second edition, original book: 1980) for example.

Thanks to the introduction of the new characterizations and the new proof of \cite{nest} and its variants,
this paper offers a \emph{homogeneous} and \emph{self-contained}\footnotemark{}
 combined presentation of the three notions, with
simpler characterizations, and short and straightforward proofs.
\footnotetext{Note that this paper does only use the most natural property of \emph{graphs}, 
that is a \emph{cardinality} argument: since edges are of maximal size two,
it is easy to see that a acyclic graph on $n$ vertices has at most $n-1$ (non-singleton) edges.
That is to say: we make no use of the other numerous connections, that are detailed in \cite{DaMoChordal},
between graph properties and hypergraph acyclicities.}

Thanks to the different characterizations of each notion, 
it becomes easier to discuss, through the study of the closure properties
of the acyclicity notions, the following questions:
\begin{itemize}
\item What makes alpha, beta and gamma acyclicities, particularly interesting --- besides their respective known applications?
\item Has every ``interesting'' hypergraph acyclicity notion been considered?
\end{itemize}
For completeness reasons, we introduce a notion of ``cycle-freedom'', that generalizes alpha acyclicity;
this notion corresponds, in \cite{Fagin83degreesof,desi} to the property of not having a \emph{pure cycle}.
Based on closure properties, we give an informal definition of ``good acyclicity notion''.
We prove that alpha, beta, and gamma acyclicity and cycle-freedom are good acyclicity notions,
and that a good acyclicity notion:
\begin{itemize}
\item either lie between gamma and beta acyclicity on one side,
\item or lie between alpha acyclicity and cycle-freedom on the other side.
\end{itemize}
%
In particular, a good acyclicity notion cannot be between beta and alpha acyclicity.
The basic fact that alpha acyclicity and beta acyclicity are the extremities of the ``gap''
makes them of particular interest.
\paragraph{Organization of this paper}
In a first section, we introduce all general-purpose definition, and in particular
a notion of cycle-freedom,
which is a naive notion (being pure-cycle free in \cite{Fagin83degreesof,desi}) more general than alpha acyclicity.
We introduce usual definitions of type (1), that we call (1a), of alpha, beta, and gamma acyclicity;
more precisely, these definitions are characterizations that make the hierarchy obvious,
for example a hypergraph is beta acyclic iff all its subset are \emph{alpha-acyclic}.
We prove alternative characterizations of type (1), called (1b), that are focused on making the considered acyclicity notion easier to establish:
for example, a hypergraph is beta acyclic iff all its subsets are \emph{cycle-free}.
These characterizations are summed up in \thref{char:type1}.
Then we introduce our new characterizations of type (2), that we call (2a),
in terms of alpha (resp. beta, gamma) leaves. 
The easy part (2a)$\Rightarrow$(1) is proved there (\thref{lem:easyway}).

In a second section, we give the technically involved results, the implications
(1)$\Rightarrow$(2a) in the three cases of alpha, beta and gamma acyclicity.
This gives the main equivalence, see \thref{char:hard}.

In a third section, we show how we can easily make the connection with other well-known
characterizations of type (2). More precisely:
\begin{description}
\item[(2b)] Characterizations in terms of non-deterministic reduction process:
GYO-reduction (\cite{Graham,YO}) for alpha acyclicity,
and the elegant DM-reduction (from \cite{82-DaMo}) in the case of gamma acyclicity.
\item[(2c)] Characterizations in terms of join tree\footnote{The reader familiar with alpha acyclicity may be surprised to find the characterization
in terms of join tree considered as a characterization ``in terms of a reduction process''.
This is discussed in the concerned section. 
In a nutshell, ``a hypergraph is GYO-reducible'' and ``a hypergraph has a join tree'' are close statements.} 
 (or of join tree with disjoint branches \cite{theseDuris,betagamma}). 
\end{description}
These characterizations are summed up in the following results:
\begin{description}
\item[alpha acyclicity] see \thref{char:alpha}, page \pageref{char:alpha}.
\item[beta acyclicity] see \thref{char:beta}, page \pageref{char:beta}.
\item[gamma acyclicity] see \thref{char:gamma}, page \pageref{char:gamma}.
\end{description}
This concludes the presentation of the acyclicity notions through their different characterizations.

Finally, the fourth last section is devoted to the study of closure properties that is presented in the introduction.
%
%
%

\section{Definitions and Properties}
\noindent
In this document, we adopt the following typographic conventions:
\begin{description}
\item[hypergraphs] or functions returning hypergraphs are written in \emph{calligraphic} font, e.g. $\Hy$, $\mathcal{G}$, $\M(\Hy)$ (minimization of $\Hy$).
\item[sets] are in \emph{uppercase} while \emph{elements} are in \emph{lowercase}, e.g. a set of vertices $S$, $\Vertexp{\Hy}$ the set of vertices of $\Hy$; 
	by contrast, vertices are always in lowercase e.g. $x$, $y$, $x_1$, etc.
	An edge is an element of some hypergraph (despite that it is also a set of vertices),
	therefore it is written in lowercase, e.g. an edge is written $e$.
	For example: ``the set of vertices $S$ is therefore an edge $e$ of the hypergraph $\Hy$, i.e. $S=e\in\Hy$.''
\item[ordered sets] (tuples) are denoted by \emph{vectors} i.e. $\vec{v}=(v_1,\ldots,v_n)$.
\end{description}
\subsection{Definitions and Type 1 Characterizations}
\noindent
First of all, 
we make a remark on the definition of a hypergraph.
\begin{remen}
Most often, people define hypergraphs as pairs $H=(V,E)$ with $V$ the set of vertices,
and $E$ the set of edges, such that:
\begin{itemize}
\item The union of edges is contained in $V$.
\item $E$ does not contain an empty edge (the empty set).
\end{itemize}
Excluding the possibility of having an empty edge is quite natural when studying acyclicity notions:
the underlying idea is that the empty edge cannot play a role in a cycle.
Nevertheless, for the same reason, we have no reason to consider the case where some vertices
are contained in no edge. 
If we exclude this case, we do not need to define $V$ separately:
it can be inferred from $E$, i.e. $V$ can be defined as the union of all the edges.
That is why we define a hypergraph as a set of non-empty sets, that is to say the hypergraph is only a set of non-empty edges.

This way, the dual (see \thref{def:dual} and \thref{rem:dual})
of a hypergraph is a hypergraph, which is not the case with the classical definition:
vertices appearing in no edge resulted in an empty edge of the dual hypergraph.

This also simplifies the notation, and make some definitions more natural:
``we say $\Hy_1$ is a subhypergraph of $\Hy_2$ when $\Hy_1\subseteq\Hy_2$''.
\end{remen}
\begin{defen}[hypergraphs and transformations] \label{def:dual}
A \emph{hypergraph} is a set of non-empty sets,
that are called its \emph{edges};
the set of \emph{vertices} of a hypergraph $\Hy$, denoted $\Vertexp{\Hy}$,
is defined as the union of all its edges.
The \emph{size} of a hypergraph $\Hy$ is defined as the sum of the cardinality of its edges i.e.
\smash {$\Sigma_{e\in\Hy}\mathrm{card}\,e$.}

A hypergraph $\Hy'$ is a \emph{subhypergraph} (or simply: \emph{subset}) of a hypergraph $\Hy$ when $\Hy'\subseteq\Hy$;
in this case we also say $\Hy'$ is obtained from $\Hy$ by \emph{removing edges}.

A hypergraph $\Hy'$ is the \emph{induced subhypergraph of $\Hy$ on a set $S\subseteq\Vertexp{\Hy}$}, denoted $\Hy'=\Hy[S]$, when
$\Hy[S]=\{e\cap S\,|\,e\in\Hy\}\smallsetminus\{\emptyset\}$;
in this case we also say $\Hy[S]$ is obtained from $\Hy$ by \emph{removing vertices} (those in $\Vertexp{\Hy}\smallsetminus S$).
For short, we write $\Hy[\smallsetminus S]$ for $\Hy[\Vertexp{\Hy}\smallsetminus S]$.

Two edges $e$ and $f$ of a hypergraph are \emph{intersecting} when we have $e\nsubseteq f$, and $f\nsubseteq e$, and $e\cap f\neq\emptyset$;
in other words, when we can find $x$, $y$ and $z$ such that \smash{$\{e,f\}\big[\{x,y,z\}\big]=\big\{\{x,y\},\{y,z\}\big\}$}.

The \emph{star of the vertex $x$ in the hypergraph $\Hy$}, denoted $\Hy(x)$, is defined as $\Hy(x)=\{e\in\Hy\,|\,x\in\Hy\}$.

The \emph{dual} of a given hypergraph $\Hy$, denoted $\mathcal{D}(\Hy)$, is defined as:
$\mathcal{D}(\Hy)=\{\Hy(x)\,|\,x\in\Vertexp{\Hy}\}$.
The \emph{normalization} of a given hypergraph $\Hy$, denoted $\mathcal{N}(\Hy)$, is defined as $\mathcal{D}(\mathcal{D}(\Hy))$.

The \emph{minimization of a hypergraph $\Hy$}, denoted $\M(\Hy)$, 
is defined as $\M(\Hy)=\{e\in\Hy\,|\,\nexists f\in \Hy\: e\subset f\}$,
that is to say the set of edges that are maximal for inclusion.

A hypergraph is a \emph{graph}, 
when all its edges have cardinality at most two, i.e. graphs may contain self-loops.
\end{defen}
Notice that applying $\M$ to a graph only removes singleton edges (usually called \emph{loops}).
\begin{figure}[t]
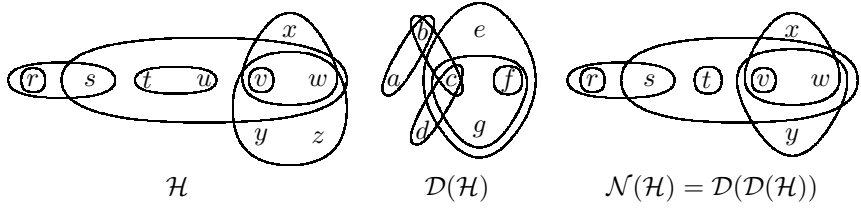

\begin{small}
\begin{center}
\begin{tabular}{ccc}
\input{ex1.tex}\vspace{3pt}
&\input{ex1dual.tex}&\input{ex2.tex}\\
$\Hy$&$\mathcal{D}(\Hy)$&$\mathcal{N}(\Hy)=\mathcal{D}(\mathcal{D}(\Hy))$
\end{tabular}
\end{center}
\end{small}
\caption{Normalization. 
Let $\Hy$ be defined as above.
Let $a=\{r\}$, $b=\{r,s\}$, $c=\{s,t,u,v,w\}$, $d=\{t,u\}$, $e=\{v\}$, $f=\{x,v,w\}$ and $g=\{v,w,y,z\}$.
Then $\Hy=\{a,b,c,d,e,f,g\}$, $\Hy(r)=\{a,b\}$,
$\Hy(s)=\{b,c\}$,
$\Hy(t)=\Hy(u)=\{d,c\}$,
$\Hy(v)=\{e,c,f,g\}$,
$\Hy(w)=\{c,f,g\}$,
$\Hy(x)=\{f\}$, and
$\Hy(y)=\Hy(z)=\{g\}$.
Therefore, the second hypergraph is the dual of the first one.\label{fig:norm}}
\vspace*{6pt}
\end{figure}
\begin{remen}	\label{rem:dual}
Informally,
the \emph{dual} of a hypergraph \smash{$\Hy=\{e_1,\ldots,e_k\}$} on vertices $\Vertexp{\Hy}=\{x_1,\ldots,x_n\}$
is obtained by ``exchanging'' the role of the vertices and of the edges, that is the hypergraph 
\smash{$\Hy'=\{y_1,\ldots,y_n\}$} (where some $y_i$ and $y_j$ may be the same edge\footnote{%
Consider the hypergraph $\{\{x,y,z\},\{y,z,t\}\}$. Its dual is isomorphic to $\{\{e\},\{e,f\},\{f\}\}$, because $y$ and $z$ are contained in the same set of edges.}) 
of set of vertices $\Vertexp{\Hy'}=\{f_1,\ldots,f_k\}$
such that for $i\leq n$ and $j\leq k$, $f_i\in y_j$ iff $x_i\in e_j$.

Technically speaking, the hypergraph we get by exchanging edges and vertices is \emph{isomorphic}
to the dual.
\end{remen}
\begin{remen}
Notice that the normalization of a hypergraph $\Hy$ is isomorphic to the hypergraph obtained from $\Hy$
by repeating the following process until it cannot be applied: ``if there are two distinct vertices $x$ and $y$ such that
$\Hy(x)=\Hy(y)$, remove the vertex $y$ from $\Hy$.
In the literature, this process is called \emph{contracting all modules}.
As an example, see Figure~\ref{fig:norm}.
\end{remen}
We give a few facts that will be used extensively.
\begin{remen} The following facts are trivial:
\label{lem:trivialfacts}
\begin{itemize}
\item We have $\M(\Hy)\subseteq\Hy$. More precisely, $\M(\Hy)$ is the set of edges of
	$\Hy$ that are maximal for set inclusion.
\item We have $\Hy[S][T]=\Hy[S\cap T]=\Hy[T][S]$ for every $T\subseteq\Vertexp{\Hy}$ and $S\subseteq\Vertexp{\Hy}$.
\item We have $\M(\M(\Hy)[S])=\M(\Hy[S])$.
\item We have $\Hy[S](x)=\Hy(x)[S]$. If $x\notin S$, then $\Hy[S](x)=\Hy(x)[S]=\emptyset$. 
\item A hypergraph $\Hy'$ can be obtained from $\Hy$ by removing vertices and edges (in any order) iff $\Hy'\subseteq\Hy[S]$ for some $S\subseteq\Vertexp{\Hy}$.
\end{itemize}
\end{remen}
\begin{defen}[hypergraph properties]
Two vertices $x$ and $y$ are \emph{neighbours}
in $\Hy$ when there is some $e\in\Hy$
such that $\{x,y\}\subseteq e$. Equivalently, $x$ and $y$ are neighbours in $\Hy$ iff $\Hy(x)\cap\Hy(y)\neq\emptyset$.
The \emph{neighbourhood} of $x$ in $\Hy$ is
the set of vertices that are neighbours to $x$, i.e. the set $\Vertexp{\Hy(x)}$\rlap{\footnotemark}.\footnotetext{%
Thanks to our definition of hypergraphs,
a vertex always belongs to its own neighbourhood; with the standard $(V,E)$ notation,
it would be necessary to deal specifically with this particular case.} 
A \emph{clique} of a hypergraph is a subset of its vertices 
whose elements are pairwise neighbours.
A hypergraph is \emph{conformal}, 
if each of its cliques is included in an edge.

A \emph{usual graph cycle} (of length $n$) is a graph $\mathcal{G}$ such that $\M(\mathcal{G})$ is isomorphic to:
\smash{$\big\{\{t_i,t_{i+1}\}\,|\,1\leq i\less n\big\}\cup\big\{\{t_n,t_1\}\big\}$}.
A tuple $\vec{t}=(t_1,\ldots,t_n)$ of $n$ pairwise distinct vertices is a \emph{cycle} of a hypergraph $\Hy$ when:
\mydisplaystyle{\begin{align*}
\M\big(\Hy[\{t_i\,|\,1\leq i\leq n\}]\big)=\big\{\{t_i,t_{i+1}\}\,|\,1\leq i\less n\big\}\cup\big\{\{t_n,t_1\}\big\}
\end{align*}}
A hypergraph $\Hy$ has a \emph{cycle} when we can find $\vec{t}$ that is a cycle of $\Hy$ or, equivalently,
when we can find $S\subseteq\Vertexp{\Hy}$ such that $\Hy[S]$ 
is a usual graph cycle (or, equivalently, $\Mp{\Hy[S]}$ is a loop-free usual graph cycle).
A hypergraph $\Hy$ is \emph{cycle-free}, 
when it has no cycle, we call \emph{cycle-freedom} this property.

A graph is acyclic when it is cycle-free.
A hypergraph $\Hy$ is \emph{alpha acyclic}, 
iff it is both conformal and cycle-free.
A hypergraph $\Hy$ is \emph{beta acyclic}, 
iff all its subsets are alpha acyclic.
A hypergraph $\Hy$ is \emph{gamma acyclic}, 
iff it is beta acyclic 
and we cannot find $x$, $y$, $z$ 
such that $\{\{x,y\},\{y,z\},\{x,y,z\}\}\subseteq \Hy[\{x,y,z\}]$.
\end{defen}
Examples of acyclic hypergraphs are given Figure~\ref{fig:exhyper}.
\begin{figure}
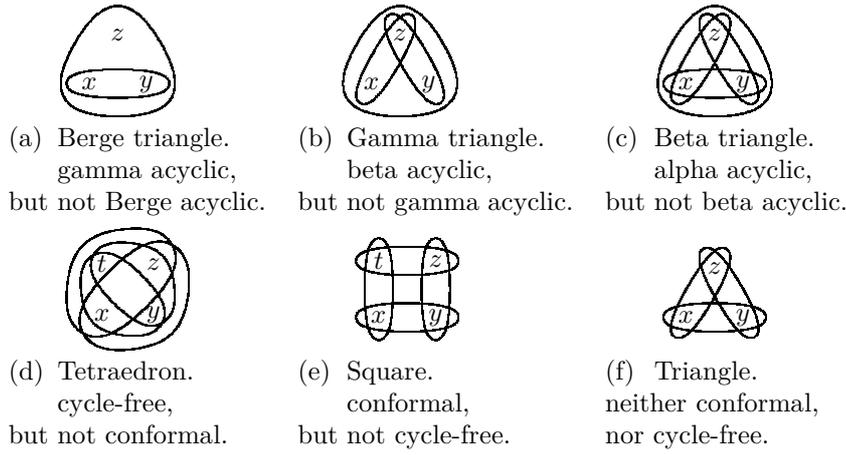

\begin{center}
\begin{small}
\begin{tabular}{@{}l@{ }ll@{ }ll@{ }l@{}}
&~\input{gamma2.tex}\vspace{3pt}&&~\input{beta2.tex}&&~\input{alpha.tex}\\ 
(a)&Berge triangle.&(b)&Gamma triangle.&(c)&Beta triangle.\\
&gamma acyclic,&&beta acyclic,&&alpha acyclic,\\
but&not Berge acyclic.&
but&not gamma acyclic.&
but&not beta acyclic.
\vspace{3pt}\\
&~\input{tetra.tex}\vspace{3pt}&&~\input{square.tex}&&~\input{triangle.tex}\\ 
(d)&Tetraedron.&(e)&Square.&(f)&Triangle.\\
&cycle-free,&&conformal,&\multicolumn{2}{l}{neither conformal,}\\
but&not conformal.&
but&not cycle-free.&
\multicolumn{2}{l}{nor cycle-free.}\\
\end{tabular}
\end{small}
\end{center}
\caption{Examples of hypergraphs. These are canonical counter-examples that we can use to prove the properties are different.
Berge acyclicity is a hypergraph acyclicity notion more restrictive than gamma acyclicity, it is defined Section~4.
 Notice that (b),(d),(e),(f) are respectively isomorphic to their own dual. 
More generally, loop-free graph cycles are self-dual, and so are graph paths that have one loop exactly, that is on one extremity of the path.\label{fig:exhyper}}
\vspace*{6pt}
\end{figure}
The following definitions will be extensively used section~4.
\begin{defen}[closure and invariance properties]\label{def:prop}
We say a property $P$ of hypergraphs is \emph{closed under} (resp. \emph{invariant w.r.t}) an operation\footnote{%
We do not require the operation to be deterministic, e.g. removing an edge. That is why
we do not write the operation as a function, and state explicitly ``$\Hy'$ is obtained from $\Hy$ by applying the operation \ldots''} $\mathcal{O}$
when, for all hypergraph $\Hy$, and all $\Hy'$ obtained from $\Hy$ by the operation $\mathcal{O}$, the following holds: $P(\Hy)\Rightarrow$ (resp. $\Leftrightarrow$) $P(\Hy')$.
\end{defen}
\begin{remen} \label{lem:stabinduc}
By the definitions, some facts are easy to derive.
These basic properties will be used \emph{extensively}.

If a hypergraph $\Hy$ is 
cycle-free (resp. conformal, alpha acyclic, beta acyclic, gamma acyclic)
then, for every $S\subseteq\Vertexp{\Hy}$, so is $\Hy[S]$.
That is to say cycle-freedom, conformity, alpha acyclicity, beta acyclicity, and gamma acyclicity are closed under taking induced hypergraphs.

A hypergraph $\Hy$ is cycle-free (resp. conformal, alpha acyclic) iff $\M(\Hy)$ is.
That is to say cycle-freedom, conformity and alpha acyclicity are invariant w.r.t. the minimization operation $\M$.

If a hypergraph is beta (resp. gamma) acyclic, then so is every subset of it.
That is to say beta acyclicity and gamma acyclicity are closed under taking a subset.
\end{remen}
We give an alternative definition of conformity, that, in some proofs, is handier to work with.
The idea is to consider a canonical clique that makes a hypergraph non-conformal.
\begin{charen}[conformity]
\label{lem:conf}
A hypergraph $\Hy$ is conformal iff
there is no $S\subseteq\Vertexp{\Hy}$, of cardinality at least $3$ such that $\M(\Hy[S])=\{S\smallsetminus\{x\}\,|\,x\in S\}$.
\begin{proofen}
If we can find $S$ such that $\M(\Hy[S])=\{S\smallsetminus\{x\}\,|\,x\in S\}$,
then obviously $\Hy$ is not conformal, by \thref{lem:stabinduc}.

If $\Hy$ is not conformal, then we can find a clique $K$ that is not included
in any edge, we assume $K$ is \emph{minimal} for this property.
Suppose $K$ has cardinality 2. The vertices in $K$ are neighbours, therefore some edge includes them both, contradiction.
So, $K$ has cardinality at least 3.
Since $K$ is minimal, any $K\smallsetminus\{x\}$ with $x\in K$ is a clique included
in an edge $e\in\Hy$, that must not include $K$, therefore $e\cap K=K\smallsetminus\{x\}$,
finally $\{K\smallsetminus\{x\}\,|\,x\in K\}\subseteq \Hy[K]$.
Since $K$ is not included in an edge of $\Hy$, $\Hy[K]$ does not contain $K$.
Finally, $\M(\Hy[K])=\{K\smallsetminus\{x\}\,|\,x\in K\}$.
\end{proofen}
\end{charen}
\begin{remen}
It is easy to see that a hypergraph $\Hy$ is conformal iff the maximal cliques of $\Hy$ are the edges of $\M(\Hy)$.
\end{remen}
We give an alternative characterization of beta acyclicity that does not refer to alpha acyclicity (and not even to conformity),
and shows how natural this notion is.
Another way of stating the following result (thanks to the last point of \thref{lem:trivialfacts}) would be:
a hypergraph $\Hy$ is beta acyclic iff \emph{we cannot get a usual graph cycle from $\Hy$ by removing
arbitrarily some edges and some vertices}.
\begin{charen}[beta acyclicity] \label{lem:beta1b} %
A hypergraph $\Hy$ is beta acyclic if and only if:
\begin{description}
\item[($\beta$1b)] Every subset of $\Hy$ is cycle-free.
\end{description}
\begin{proofen}
If some subset $\Hy'$ of a hypergraph $\Hy$ is not cycle-free, then $\Hy'$
is not alpha acyclic, therefore $\Hy$ is not beta acyclic.

Conversely, if $\Hy$ is not beta acyclic, then we can find a subset $\Hy'\subseteq\Hy$
that is not alpha acyclic.
Suppose $\Hy'$ is cycle-free, then $\Hy'$ is not conformal.

Then, for some $S$ with at least three elements, $\Hy''=\M(\Hy'[S])=\{S\smallsetminus\{x\}\,|\,x\in S\}$.
Assume $S$ has three elements. Then $\Hy''$ is not cycle-free, therefore $\Hy'$ is not, contradiction.

Take $x,y,z$ in $S$. We have $\Hy''[\{x,y,z\}]=\{\{x,y\},\{y,z\},\{x,z\},\{x,y,z\}\}$,
The hypergraph $\Hy_2'\subset\Hy'$, defined as the set of edges of $\Hy'$
that do not include $\{x,y,z\}$ satisfies $\Hy_2'[\{x,y,z\}]=\{\{x,y\},\{y,z\},\{x,z\}\}$,
therefore $\Hy_2'\subset\Hy$ is not cycle-free, by \thref{lem:stabinduc}, a contradiction.
\end{proofen}
\end{charen}
We generalize a result of Duris \cite{theseDuris,betagamma} that states that a hypergraph $\Hy$ is gamma acyclic iff it is \emph{alpha} acyclic and we cannot find $x,y,z$ such that $\{\{x,y\},\{x,z\},\{x,y,z\}\}\subseteq \Hy[\{x,y,z\}]$. This ``weaker'' characterization allows to simplify some proofs.
\begin{charen}[gamma acyclicity]
\label{char:gammaacy}
A hypergraph $\Hy$ is gamma acyclic if and only if:
\begin{description}
\item[($\gamma$1b)] $\Hy$ is cycle-free and we cannot find $x,y,z$ such that:
\mydisplaystyle{\begin{align*}\{\{x,y\},\{x,z\},\{x,y,z\}\}\subseteq\Hy[\{x,y,z\}]\end{align*}}
\end{description}
\begin{proofen}
%
%
We only prove the non-trivial part of the result, that is to say: 
if a hypergraph $\Hy$ is cycle-free but not gamma acyclic, then we can find $x,y,z$ such that 
$\{\{x,y\},\{x,z\},\{x,y,z\}\}\subseteq \Hy[\{x,y,z\}]$.
Let $\Hy$ be a cycle-free but not gamma acyclic hypergraph.
If $\Hy$ is beta-acyclic, the result is obvious.
It remains to prove that, if $\Hy$ is a cycle-free but not beta acyclic hypergraph
then we can find $x,y,z$ such that 
$\{\{x,y\},\{x,z\},\{x,y,z\}\}\subseteq \Hy[\{x,y,z\}]$.

Let $S\subseteq\Vertexp{\Hy}$ such that $\Hy[S]$ is not beta acyclic, we assume $S$ minimal for this property.
Let $\Hy'\subseteq\Hy[S]$ such that $\Hy'$ is not cycle-free, we assume $\Hy'$ minimal.
We can find $S'$ such that $\Hy'[S']$ is a usual graph cycle.
Assuming $S'\neq S$ contradicts that $S$ is minimal, therefore $\Hy'=\Hy'[S']$ itself is a usual graph cycle.
Some edge $e\in\Hy[S]\smallsetminus\Hy'$, has cardinality at least two, or else $\Hy[S]$ would not be cycle-free.
We assume w.l.o.g. that $e$ contains at least the vertex $x_1$, and 
that $\Hy'$ is in the form $\{\{x_i,x_{i+1}\,|\,1\leq i< n\}\cup\{x_1,x_n\}$.

Assume $e\neq S$.
Let $x_i$ be the vertex of smallest subscript that is not contained in $e$.
Assume $i>2$.
In this case, 
$(\Hy'\cup\{e\})[\{x_{i-1},\ldots,x_n\}]$ is a graph cycle,
therefore $\Hy[\{x_{i-1},\ldots,x_n\}]$ is not beta acyclic, which contradicts that $S$ is the minimal.
Therefore, $x_1$ in $e$ but $x_2$ not in $e$.
Now, since $e$ has cardinality $2$ at least, it contains another vertex.
Let $x_i$ be the smallest vertex that is in $e$. Since $e$ is not an edge of the cycle,
$i\neq n$. But then, $(\Hy'\cup\{e\})[\{x_{1},\ldots,x_i\}]$ is a graph cycle,
and once more we get a contradiction with that $S$ is minimal.

We have proved $e=S$.
Notice that 
as a consequence, $(\Hy'\cup\{e\})[\{x_1,x_2,x_3\}]$ contains $\{x_1,x_2\}$, $\{x_2,x_3\}$, and $\{x_1,x_2,x_3\}$.
We have proved $\Hy[S][\{x_1,x_2,x_3\}]=\Hy[\{x_1,x_2,x_3\}]$ contains $\{x_1,x_2\}$, $\{x_2,x_3\}$, and $\{x_1,x_2,x_3\}$.
\end{proofen}
\end{charen}
\begin{remen}
\label{rem:acygraph}
On \emph{graphs}, acyclicity, cycle-freedom, alpha acyclicity, beta acyclicity, and gamma acyclicity are equivalent.
\begin{proofen}
Considering the definitions, 
we only have to prove that a non gamma acyclic graph is not acyclic.
By the definitions, gamma acyclicity and beta acyclicity are equivalent when considering graphs.
Let $\mathcal{G}$ be a non beta acyclic graph.
Some subset of $\mathcal{G}$ has a cycle. It is easy to check this cycle is a cycle of $\mathcal{G}$.
\end{proofen}
\end{remen}
\begin{remen}
\begin{enumerate}
\item There are alpha acyclic hypergraphs whose subsets are not all alpha acyclic: e.g. $\{\{x,y\},\{y,z\},\{x,z\},\{x,y,z\}\}$.
\item A subset of a beta (resp. gamma) acyclic hypergraph is a beta (resp. gamma) acyclic hypergraph, by definition.
\item A subset of an acyclic graph is an acyclic graph, since the subset of a graph is a graph and by the previous point.
\end{enumerate}
\end{remen}
Let us sum up the results proved so far.
We have proved the equivalence of different characterizations of type (1) for each degree of acyclicity.
\begin{charen}[first type characterizations]
\label{char:type1}
Let $\Hy$ be a hypergraph.

The following are equivalent:
\begin{description}
\item[($\alpha$1a)] The hypergraph $\Hy$ is alpha acyclic, i.e. conformal and cycle-free. (standard)
\item[($\alpha$1b)] We \emph{cannot} find $S\subseteq\Vertexp{\Hy}$ such that $\M(\Hy[S])$ is either $\{S\smallsetminus\{x\}\,|\,x\in S\}$ or a usual graph cycle.
\end{description}
The following are equivalent:
\begin{description}
\item[($\beta$1a)] The hypergraph $\Hy$ is beta acyclic, i.e. every subset of $\Hy$ is alpha acyclic. (standard)
\item[($\beta$1b)] Every subset of $\Hy$ is cycle-free.
\item[($\beta$1b)] (reformulated) We cannot get a usual graph cycle from $\Hy$ by removing vertices and/or edges.
\end{description}
The following are equivalent:
\begin{description}
\item[($\gamma$1a)] The hypergraph $\Hy$ is gamma acyclic, i.e. $\Hy$ is beta acyclic and we cannot find $x,y,z$ such that $\{\{x,y\},\{x,z\},\{x,y,z\}\}\subseteq \Hy[\{x,y,z\}]$. (our definition)
\item[($\gamma$1b)] The hypergraph $\Hy$ is cycle-free and we cannot find $x,y,z$ such that $\{\{x,y\},\{x,z\},\{x,y,z\}\}\subseteq \Hy[\{x,y,z\}]$. 
(Definition~3 of gamma acyclicity in \cite{Fagin83degreesof})
\end{description}
\begin{proofen}
The equivalence ($\alpha$1a)$\Leftrightarrow$($\alpha$1b) is proved by \thref{lem:conf}.

The equivalence ($\beta$1a)$\Leftrightarrow$($\beta$1b) is proved by \thref{lem:beta1b}.

The equivalence ($\gamma$1a)$\Leftrightarrow$($\gamma$1b) is proved by \thref{char:gammaacy}.
%
\end{proofen}
\end{charen}
Notice that the equivalence ($\gamma$1a)$\Leftrightarrow$($\gamma$1b) 
shows that the exclusion of the ``gamma triangle'' (see Figure~\ref{fig:exhyper}), 
i.e. the pattern $\{\{x,y\},\{x,z\},\{x,y,z\}\}$, is quite a constraining statement.
\subsection{Leaves and Elimination Orders}
\noindent
We introduce the notion of alpha (resp. beta, gamma) \emph{leaves},
and the associated notions of alpha (resp. beta, gamma) \emph{elimination orders},
that are characterizations of type (2), i.e. in terms of reducibility to the
empty hypergraph through a certain reduction process.
We prove the existence of an alpha (resp. beta, gamma) elimination order implies alpha (resp. beta, gamma)
acyclicity.
\begin{defen}[leaves, eliminations orders]
A vertex $x$ of $\Hy$ is an \emph{alpha leaf}, 
if $\M(\Hy(x))$ has a single edge,
that is to say $\Hy(x)$ has a maximal element for inclusion.

A vertex $x$ of $\Hy$ is a \emph{beta leaf}, 
if for all $e$, $f$ in $\Hy( x)$, $e\subseteq f$ or $f\subseteq e$.

A vertex $x$ of $\Hy$ is a \emph{gamma leaf}, 
if:
\begin{itemize}
\item The vertex $x$ is a beta leaf of $\Hy$, we call $e_x$ the maximal edge holding $x$ in $\Hy$.
\item Every neighbour of $x$ in $\Hy\smallsetminus\{e_x\}$ is a beta leaf.
\end{itemize}
 
We say $\vec{v}=(v_1,\ldots,v_n)$ is a \emph{alpha} (resp. \emph{beta}, \emph{gamma}) \emph{elimination order} of a hypergraph $\Hy$
when $v_n$ is an alpha leaf (resp. beta leaf, gamma leaf) 
and $(v_1,\ldots,v_{n-1})$ is an alpha (resp. beta, gamma) elimination order of $\Hy[\smallsetminus\{x_n\}]$.
The empty tuple is an alpha (resp. beta, gamma) elimination order of the empty hypergraph.
\end{defen}
\begin{remen}
A good way to understand what an alpha (resp. beta, gamma) leaf is, is to understand what it is not.
The following facts, that will be extensively used, are easily derived from the definitions:
\begin{itemize}
\item A vertex $x$ is \emph{not} an alpha leaf of a hypergraph $\Hy$ iff there are two vertices $y$ and $z$ such that $\M(\Hy[\{x,y,z\}])=\{\{x,y\},\{x,z\}\}$;
	in other words, $\M(\Hy(x))=(\M(\Hy))(x)$ has two intersecting edges.
\item A vertex $x$ is \emph{not} a beta leaf of a hypergraph $\Hy$ iff there are two vertices $y$ and $z$ such that $\{\{x,y\},\{x,z\}\}\subseteq\Hy[\{x,y,z\}])$;
	in other words, $\Hy(x)$ has two intersecting edges.
\item A vertex $x$ is \emph{not} a gamma leaf of a hypergraph $\Hy$ iff either it is not a beta leaf or there are two vertices $y$ and $z$ such that 
	$\{\{x,y\},\{y,z\},\{x,y,z\}\}\subseteq\Hy[\{x,y,z\}])$.
\end{itemize}
Additionally, observe that a beta leaf $x$ of $\Hy$ is an alpha leaf of every $\Hy'\subseteq\Hy$ such that $x\in\Vertexp{\Hy'}$.
This is to be put in relation with \thref{lem:beta1b}, i.e. the characterization ($\beta$1b) of beta acyclicity, that states: 
a hypergraph is beta acyclic iff all its subset are alpha acyclic.
\end{remen}
%
%
The following is easy but long to prove formally.
\begin{lemen}
\label{lem:eqleaf}
Let $\Hy$ be a hypergraph, having an alpha (resp. beta, gamma) leaf denoted $x$.
Then $\Hy$ is alpha (resp. beta, gamma) acyclic iff $\Hy[\smallsetminus\{x\}]$ is.
\begin{proofen}
Let $\Hy$ be a hypergraph, having an alpha (resp. beta, gamma) leaf we call $x$.
We already know, by \thref{lem:stabinduc}, that if $\Hy$ is alpha (resp. beta, gamma) acyclic,
then so is $\Hy[\smallsetminus\{x\}]$. 

Assume that $x$ is an alpha leaf, and that $\Hy[\smallsetminus\{x\}]$ is alpha acyclic.
Assume $\Hy$ is not conformal. 
Then, by \thref{lem:conf}, there is some $S\subseteq\Vertexp{\Hy}$, of cardinality at least 3, such that $\M(\Hy[S])=\{S\smallsetminus\{v\}\,|\,v\in S\}$.
Assume $x\notin S$, then $\Hy[S]=\Hy[\smallsetminus\{x\}][S]=\Hy'[S]$ is not alpha acyclic,
therefore $\Hy'$ is not, a contradiction.
As a consequence, $x\in S$. Since $x$ is an alpha leaf of $\Hy$, 
there is an edge of $\Hy$ that includes every neighbour of $x$ (including $x$ itself);
this edge therefore includes $S$, therefore $S\in\Hy[S]$, a contradiction.
We have proved $\Hy$ is conformal.
With the same reasoning, we prove $\Hy$ is cycle-free.
Assume this is not the case: for some $S=\{x_1,\ldots,x_k\}$,
$\M(\Hy[\{x_1,\ldots,x_k\}])=\{\{x_i,x_{i+1}\}\,|\,1\leq i<k\}\cup\{\{x_1,x_k\}\}$.
If $x\notin S$, we have the same contradiction as before; hence $x\in S$,
let us assume w.l.o.g $x$ is $x_1$.
Then, since $x$ is an alpha leaf of $\Hy$, some edge includes all his neighbours, in particular,
it includes both $x_2$ and $x_k$, which contradicts the fact that they are not neighbours in $\M(\Hy[S])$.
We have proved $\Hy$ is alpha acyclic.

Assume that $x$ is a beta leaf, and that $\Hy[\smallsetminus\{x\}]$ is beta acyclic and that $\Hy$ is not.
In this case, thanks to characterization ($\beta$1b) of beta acyclicity,
we can find $\Hy'\subseteq\Hy$ and $S=\{x_1,\ldots,x_k\}\subseteq\Vertexp{\Hy'}$
such that: $\M(\Hy[\{x_1,\ldots,x_k\}])=\{\{x_i,x_{i+1}\}\,|\,1\leq i<k\}\cup\{\{x_1,x_k\}\}$.
If $x\notin S$, then $\Hy'[\smallsetminus\{x\}]\subseteq\Hy[\smallsetminus\{x\}]$ 
is not beta acyclic, a contradiction.
Assume w.l.o.g. $x$ is $x_1$. Then there are two incomparable edges of $\Hy'$, which are edges of $\Hy$,
that include $x$, therefore $x$ is not a beta leaf of $\Hy$, a contradiction.
We have proved $\Hy$ is beta acyclic.

Assume that $x$ is a gamma leaf, and that $\Hy[\smallsetminus\{x\}]$ is gamma acyclic.
By the previous point, and since a gamma leaf is a beta leaf, we already know that $\Hy$ is beta acyclic.
We only have to prove that there is no vertices $s,t,u$ such that 
$\{\{s,t,u\},\{s,t\},\{t,u\}\}\subseteq\Hy[\{s,t,u\}]$.
Assume we can find such $s,t,u$. 
If $x$ is none of them, then $\Hy'$ is not gamma acyclic, a contradiction.
If $x$ is $t$, then $x$ is not a beta leaf, a contradiction.
By symmetry, assuming $x$ is $s$ and assuming $x$ is $u$ is the same, we assume $x$ is $s$.
In $\Hy$, the maximal edge containing $x$, called $e_x$, includes $x$, $t$ and $u$,
there is also an edge $e$ that contains $x$ and $t$ but not $u$, and an edge $f$ that contains $t$ and $u$ but not $x$.
In $\Hy\smallsetminus\{e_x\}$, the vertex $t$ is a neighbour of $x$, but is contained in the two incomparable
edges $e$ and $f$, i.e. it is not a beta leaf; this contradicts the fact that $x$ is a gamma leaf.
\end{proofen}
\end{lemen}
\begin{coren}
\label{lem:easyway}
If a hypergraph $\Hy$ has an alpha (resp. beta, gamma) elimination order,
then it is alpha (resp. beta, gamma) acyclic.
\begin{proofen}
By the previous lemma, 
if a hypergraph $\Hy$ has an alpha (resp. beta, gamma) elimination order,
then it is alpha (resp. beta, gamma) acyclic iff the empty hypergraph is, which is the case.
\end{proofen}
\end{coren}
\section{Main Results}
\noindent
In a first subsection, 
we prove the ``hard part'' of the equivalence between ``being alpha (resp. beta, gamma) acyclic'' and ``having an alpha (resp. beta, gamma) elimination order''.
This hard part consists in proving that an alpha (resp. beta, gamma) acyclic hypergraph has an alpha (resp. beta, gamma) leaf.
Then we show how this leads straightforwardly to the equivalence.
\subsection{Acyclicity Implies Having Leaves}
\noindent
The following technical lemma is folklore (it is an obvious corollary of a well-known fact), 
we prove it for the sake of self-containment, and in order to get absolutely sure that there is nothing ``complicated'' or non-straightforward hidden inside.
\begin{lemen}
\label{lem:folkcycle}
A graph where each vertex has exactly two neighbours has a cycle.
\begin{proofen}
We prove it by induction on the cardinality of $\Vertexp{\mathcal{G}}$.
Assume this is the case for every graph with $n$ vertices.
Take $\mathcal{G}$ with $n+1$ vertices.

Take any vertex $x$. It has exactly two neighbours $y$ and $z$.
If $y$ and $z$  are neighbours,
then $(x,y,z)$ is a cycle of $\mathcal{G}$. From now on, we assume the contrary.
Consider
\smash{$\mathcal{G}'=\mathcal{G}\smallsetminus\big\{\{x,y\},\{x,z\}\big\}\uplus\big\{\{y,z\}\big\}$}.
In the graph $\mathcal{G}'$, each vertex has exactly two neighbours; furthermore,
this latter graph has $n$ vertices, hence, by induction, it contains a cycle $\vec{c}$.
Either both $y$ and $z$ belong to $\vec{c}$ or not.
\begin{itemize}
\item If $y$ and $z$ both belong to $\vec{c}$, then $\vec{c}$ is in the form 
$(c_1,\ldots,c_i,y,z,c_{i+3},\ldots)$.
Then
$(c_1,\ldots,c_i,y,x,z,c_{i+3},\ldots)$ is a cycle of $\mathcal{G}$.
\item In the other case, the cycle $\vec{c}$ of $\mathcal{G}'$ is also a cycle of $\mathcal{G}$.
\end{itemize}
In all cases, $\mathcal{G}$ has a cycle, QED.
\end{proofen}
\end{lemen}
The basic idea is to prove that an alpha acyclic hypergraph has an alpha leaf by induction.
However, formulated this way, this property is very hard to prove.
We strengthen the induction hypothesis in the same way Brouwer and Kolen \cite{nest}
did in order to prove a beta acyclic hypergraph has a beta leaf.
They prove, by induction, that a beta acyclic hypergraph with at least two vertices has two different leaves.

This works, but, still, the proof is long and not easy.
In a nutshell, this is because the definition of ``different'' does not work well:
if we have two different leaves $x$ and $y$ that are included in exactly the same edges, 
they are different, but behave exactly as if they were the same.
Thus we do not fully take profit from the fact they are different. 

Rather than ``two leaves'', we want ``two \emph{really different} leaves''.
An efficient way of getting two leaves that are \emph{really different} is to require that they are not neighbours.
Now, the problem is: what if all the vertices are pairwise neighbours, i.e. the set of all edges is a clique?
All the acyclicity notions imply conformity, therefore this setting only occurs when there is an edge that contains all vertices.
We can therefore consider only hypergraphs where no edge contains all vertices, and deal specifically with this case.

The following does it, and is trivial.
\begin{lemennp}
\label{lem:trivialalpha}
A hypergraph $\Hy$ with a full edge, i.e. such that $\Vertexp{\Hy}\in\Hy$,
is alpha acyclic and has only alpha leaves.
\end{lemennp}
A nice side-effect of the strengthening of the induction hypothesis, besides the simplifications it permits, is that 
the result that is proved implies a nice property stated by \thref{thm:sacred}.
\begin{thmen}
\label{thm:alphaleaves}
A alpha acyclic non-empty hypergraph $\Hy$ such that $\Vertexp{\Hy}\notin\Hy$
has two alpha leaves that are not neighbours.
\begin{proofennoend}
This is proved by total induction on the hypergraph \emph{size}. 
The property is obvious for the hypergraph of size one: $\{\{x\}\}$.
We assume the result holds for every hypergraph of size $k\less n$.
Let $\Hy$ be an alpha acyclic hypergraph of size $n$ such that $\Vertexp{\Hy}\notin\Hy$.
First of all, we simplify the situation.

If $\Hy\neq\Mp{\Hy}$, then $\Mp{\Hy}$ is an alpha acyclic hypergraph smaller than $\Hy$, without a full edge, which has therefore two alpha leaves by induction.
These leaves are leaves of $\Hy$.
From now on we assume \smash{$\Mp{\Hy}=\Hy$}.
Now suppose that $\Hy$ has two vertices. By previous point, it is isomorphic to $\{\{x\},\{y\}\}$, therefore it satisfies the property.
From now on, we assume w.l.o.g. that $\Hy$ has at least \emph{three} vertices.

%
%
%
We will first prove $\Hy$ has an alpha leaf, and then we will prove $\Hy$  has two alpha leaves that are not neighbours.
\end{proofennoend}
\begin{proofennoend}[\smash{$\Hy$} has an alpha leaf]
Assume $\Hy$ has no alpha leaf, consider the directed graph:
\smash{$\mathcal{G}=\left\{s\rightarrow t\,\left|\,t\text{ is an alpha leaf of }\Hy\big[\smallsetminus\{s\}\big]\right\}\right.$}.
We say $x$ is a predecessor of $y$ (resp. $y$ is a successor of $x$) when $(x\rightarrow y)\in\mathcal{G}$.
Let us prove that every vertex has two predecessors and two successors.

Assume $x$ has $n$ predecessors \smash{$x_1,\ldots,x_n$} in $\mathcal{G}$.
Then $x$ is an alpha leaf in
 \smash{$\Hy[\smallsetminus\{x_i\}]$}, therefore there is an edge \smash{$e_{i\!}$}
that includes all the neighbours of $x$ in \smash{$\Hy[\smallsetminus\{x_i\}]$}, 
i.e. all the neighbours of $x$ in $\Hy$ except \smash{$x_i$}.

We call $E$ the set of neighbours of $x$ in $\Hy$.
For all $i$, we define \smash{$E_i=E\smallsetminus\{x_i\}$}.
Assume that, for a given $i$, \smash{$E_i\cup\{x_i\}$} is in $\Hy$.
Then \smash{$E_i\cup\{x_i\}$} also includes all the neighbours of $x$ in $\Hy$, therefore $x$ is a leaf of $\Hy$, contradictory.
We have proved that every \smash{$E_{i\!}=e_{i\!}$} is in $\Hy$, and that $E\notin\Hy$.

If $n\geq 3$, $E$ is a clique: \smash{$e_1$} is a clique and we have to prove that all other vertices are neighbours of \smash{$x_1$}. 
The edge \smash{$e_2$} proves that all the vertices of \smash{$e\smallsetminus\{x_2\}$} are neighbours of \smash{$x_1$}. 
The edge \smash{$e_3$} proves that \smash{$x_1$} and \smash{$x_2$} are neighbours. 
The set $E$ is a clique that is not included in an edge, therefore $\Hy$ is not conformal, hence not alpha acyclic, a contradiction.

We have proved that each vertex of $\mathcal{G}$ has at most two predecessors.
We prove it has at least two successors:
for every $x$, $\Hy'=\Hy[\smallsetminus\{x\}]$ either is such that $\Vertexp{\Hy'}\in\Hy'$ or not.
In the first case, by \thref{lem:trivialalpha}, all vertices (they are at least two, see preliminary) of $\Hy'$ are alpha leaves.
In the other case, the induction hypothesis gives the result.
Consequently, every vertex has exactly two predecessors and two successors.
Now we prove every vertex is in exactly two edges of $\Hy$, that are of cardinality two.

Let $x$ be a vertex, and let $y$ and $z$ be its two predecessors.
The vertex $x$ is in the two edges \smash{$e_{y\!}$} and \smash{$e_{z\!}$}, the two edges including respectively $\{x,y\}$ and $\{y,z\}$ of maximal cardinality.
We have \smash{$e_z\smallsetminus\{x\}\subseteq e_y$} and \smash{$e_y\smallsetminus\{x\}\subseteq e_z$}.
Therefore we have \smash{$e_y=\{y\}\cup e$} and \smash{$e_z=\{z\}\cup e$}, where \smash{$e=e_y\cap e_z$}.
If there is another edge $f$ containing $x$, necessarily incomparable with \smash{$e_{y\!}$} and \smash{$e_{z\!}$},
then \smash{$e=e_y\smallsetminus\{y\}\subset f$} and similarly for $z$. 
Then, 
\smash{$e_y=e\cup\{y\}$}, 
\smash{$e_z=e\cup\{z\}$}, and
\smash{$f=e\cup\{t,\ldots\}$}  are pairwise incomparable, and \smash{$e_{z\!}$} and $f$ are in \smash{$\Hy[\smallsetminus\{y\}]$}
and are also incomparable, therefore $x$ is not an alpha leaf of \smash{$\Hy[\smallsetminus\{y\}]$}, a contradiction.

Therefore any vertex $x$ belongs to exactly two edges, each of which holds a predecessor of $x$.
Assume another vertex $t$ is in the intersection of these two edges,
then $x$ and $t$ satisfy \smash{$\Hy(x)=\Hy(t)$}.
Let $u$ be a successor of $x$.
Let $e_u$ be the maximal edge holding $u$ in $\Hy[\smallsetminus\{x\}]$,
it contains every neighbour of $u$ in $\Hy$ except $x$, therefore it contains $t$;
this contradicts \smash{$\Hy(x)=\Hy(t)$}.
Therefore every vertex belongs to exactly the two edges $\{x,y\}$ and $\{x,z\}$, 
where $y$ and $z$ are the predecessors of $x$.

Finally, we have proved that $\Hy$ is a \emph{graph} where each vertex has exactly two neighbours;
by \thref{lem:folkcycle}, this graph has a cycle, therefore it cannot be alpha acyclic, a contradiction.
\end{proofennoend}
\begin{proofen}[\smash{$\Hy$} has two alpha leaves that are not neighbours]
We know by previous point that $\Hy$ has an alpha leaf $x$. 
Let \smash{$e_{x\!}$} be the maximal edge of $\Hy$ containing $x$,
and let $\Hy'=\Hy[\smallsetminus\{x\}]$. 
We distinguish two cases:
\begin{description}
\item[$\Vertexp{\Hy'}\in\Hy'$] 
	Since \smash{$e_{x\!}\neq\Vertexp{\Hy}$}, $\Vertexp{\Hy'}\smallsetminus e_{x\!}\neq \emptyset$. 
	Let us take $y\in\Vertexp{\Hy'}\smallsetminus e_x$. 
	In this case, $y$ is an alpha leaf of $\Hy$ that is not neighbour to $x$ in $\Hy$, the expected result.
\item[$\Vertexp{\Hy'}\notin\Hy'$] In this case, by induction, $\Hy'$ has two leaves $y$ and $z$ that are not neighbours in $\Hy'$.
Since $y$ and $z$ are not neighbours in $\Hy'$, and since \smash{$e_x\smallsetminus\{x\}\neq\Vertexp{\Hy'}$},
one of them is not in $e_x\smallsetminus\{x\}$, say $y$. 
This vertex $y$ is not a neighbour of $x$ in $\Hy$,
therefore $\Hy(y)=\Hy'(y)$,
hence $y$ is an alpha leaf of $\Hy$ that is not a neighbour of~$x$.
\end{description}
In both cases, we have found a vertex $y$ that is an alpha leaf of $\Hy$ that is not neighbour to $x$.
\end{proofen}
\end{thmen}
The following proof both gives a stronger result and has a simpler\footnote{For example, see the trick consisting in introducing the set called $I$
in the proof in \cite{nest}. Our proof is straightforward, and much shorter.} proof than the one in \cite{nest}.
Notice that this new proof becomes simple enough to use exactly the same argument to prove that (i) there is a beta leaf,
and to prove that (ii) there is another one.

But first of all, as we did before, we state a trivial lemma.
\begin{lemennp}
\label{lem:trivialbeta}%
A beta leaf of $\Hy\smallsetminus\{\Vertexp{\Hy}\}$ is a beta leaf of $\Hy$.
\end{lemennp}
\begin{thmen}%
\label{thm:betaleaves}%
A beta acyclic hypergraph $\Hy$ with at least two vertices has two beta leaves that are not neighbours in $\Hy\smallsetminus\{\Vertexp{\Hy}\}$.
\begin{proofen}
We proceed by induction on the size of the hypergraph.
Assume the result holds for every hypergraph of size $n$ or less.
Consider $\Hy$ a beta acyclic hypergraph of size $n+1$ with two vertices or more.
If $\Hy=\{\Vertexp{\Hy}\}$, the result is trivial; from now on we assume this is not the case.

Assume $\Vertexp{\Hy}\in\Hy$, let $\Hy'=\Hy\smallsetminus\{\Vertexp{\Hy}\}$, 
which is a beta acyclic hypergraph with at least one vertex.
If $\Vertexp{\Hy}=\Vertexp{\Hy'}$, then by induction it has two beta leaves that are not neighbours in
$\Hy'\smallsetminus\{\Vertexp{\Hy'}\}=\Hy\smallsetminus\{\Vertexp{\Hy}\}$, these beta leaves are beta leaves of $\Hy$ by \thref{lem:trivialbeta},
which concludes the case.
Otherwise, (if $\Vertexp{\Hy'}\neq\Vertexp{\Hy}$), 
a vertex $x\in\Vertexp{\Hy}\smallsetminus\Vertexp{\Hy'}$ 
only belongs to the edge  $\Vertexp{\Hy}\in\Hy$, hence it is a beta leaf.

Besides, $\Hy'$ has at least one vertex; if it has only one, it is obviously a beta leaf;
otherwise, by induction, $\Hy'$ has beta leaves.
In both cases, $\Hy'$ has a beta leaf $y$ that is also a beta leaf of $\Hy$.

We have proved that $\Hy$ has two beta leaves $x$ and $y$ that are not neighbours in $\Hy\smallsetminus\Vertexp{\Hy}$.
From now on, we assume $\Vertexp{\Hy}\notin\Hy$.
Furthermore, if $\Hy$ only has two vertices, then $\Hy$ is isomorphic to $\{\{x\},\{y\}\}$, which satisfies the result.
From now on, let us suppose that $\Hy$ has three vertices or more.
This concludes the preliminary of the proof.

We prove the following fact:
\begin{description}
\item[fact] If $\Hy$ has a \emph{alpha}-leaf $x$, i.e. $\Hy( x)$ has a maximal element \smash{$e_{x\!}$} for set inclusion,
then $\Hy$ has a \emph{beta}-leaf that is not neighbour of $x$ in $\Hy\smallsetminus\{\Vertexp{\Hy}\}$.
\end{description}
Let $\Hy'=\Hy[\smallsetminus\{x\}]$. 
Since $e_x\neq\Vertexp{\Hy}$ (by the previous point), we have $e_x\smallsetminus\{x\}\neq\Vertexp{\Hy[\smallsetminus\{x\}]}$.
The hypergraph $\Hy'$ is beta acyclic and therefore, by induction, has two beta leaves $y$ and $z$
that are not neighbours in $\Hy'\smallsetminus\Vertexp{\Hy'}$, 
in particular in the edge $e_x\smallsetminus\{x\}$ in $\Hy'$. 
One of them, say $y$, is not in this latter edge,
and \emph{a fortiori} not in the edge \smash{$e_{x\!}$} of $\Hy$.
By the definition of \smash{$e_{x\!}$}, the vertex $y$ is not a neighbour of $x$ in $\Hy$. 
Consequently, \smash{$\Hy( y)=\Hy'( y)$},
so the vertex $y$ is a beta leaf of $\Hy$ and is not neighbour of $x$ in $\Hy$, 
and therefore not in $\Hy\smallsetminus\{\Vertexp{\Hy}\}$.

We have proved the fact, now we use it (twice!).
The hypergraph $\Hy$ is beta acyclic, and therefore alpha acyclic,
therefore it has an alpha leaf $x$. By the fact proved above,
we can find a beta leaf $y$.
This beta leaf is also trivially an alpha leaf,
so we can apply the fact once more to $y$.
We get a vertex $z$ that is also a beta leaf, such that $y$ and $z$ are not neighbours in $\Hy$,
which concludes the proof.
\end{proofen}
\end{thmen}
\begin{lemen}
\label{lem:trivialgamma}
A hypergraph with a full edge is gamma acyclic iff it is a set of pairwise non-intersecting edges.
In particular, it has only gamma leaves.
\begin{proofen}[sketch]
We only prove that a gamma acyclic hypergraph with a full edge is a set of pairwise non-intersecting edges.
Let $\Hy$ be a gamma acyclic hypergraph such that $\Vertexp{\Hy}\in\Hy$.
Assume $\Hy$ contains two intersecting edges $e$ and $f$, i.e. $e$ and $f$ are incomparable for set inclusion and such that $e\cap f\neq\emptyset$.
Let $x\in e\cap f$, $y\in e\smallsetminus f$ and $z\in f\smallsetminus e$.
The hypergraph $\Hy[\{x,y,z\}]$ contains $\{x,y,z\}$, $\{x,y\}$ and $\{x,z\}$,
therefore $\Hy$ is not gamma acyclic, a contradiction.
We have proved any two edges are non-intersecting.
\end{proofen}
\end{lemen}
\begin{thmen}
\label{thm:gammaleaves}
A gamma acyclic hypergraph $\Hy$ with at least two vertices has two gamma leaves that are not neighbours in $\Hy\smallsetminus\{\Vertexp{\Hy}\}$.
\begin{proofen}
We proceed by induction on the size of the hypergraph.
Assume the result holds for every hypergraph of size $n$ or less.
Consider $\Hy$ a gamma acyclic hypergraph of size $n+1$ with two vertices or more.
If $\Hy=\{\Vertexp{\Hy}\}$, or if $\Hy$ has two vertices, the result is trivial; from now on we assume this is not the case.

Assume $\Vertexp{\Hy}\in\Hy$. By previous lemma, \thref{lem:trivialgamma}, $\Hy$ is a set of pairwise non-intersecting edge, and has only gamma leaves.
Since $\Hy\neq\{\Vertexp{\Hy}\}$, the hypergraph $\Hy\smallsetminus\{\Vertexp{\Hy}\}$ 
is not empty, its has vertices. 
Consider any of them $x$. 
It is a beta leaf of $\Hy\smallsetminus\{\Vertexp{\Hy}\}$, therefore some edge $e$ includes
all its neighbourhood. Take $y\in \Vertexp{\Hy}\smallsetminus e$.
We have $x$ and $y$ two gamma leaves of $\Hy$ that are not neighbours in $\Hy\smallsetminus\{\Vertexp{\Hy}\}$.
From now on, we assume $\Vertexp{\Hy}\notin\Hy$.

We prove the following fact: 
\begin{description}
\item[fact] If $\Hy$ has a \emph{beta} leaf $x$, in particular $\Hy( x)$ has a maximal element \smash{$e_{x\!}$} for set inclusion,
then $\Hy$ has a \emph{gamma} leaf that is not a neighbour of $x$ in $\Hy\smallsetminus\{\Vertexp{\Hy}\}$.
\end{description}
Let $\Hy'=\Hy[\smallsetminus\{x\}]$. 
Assume $\Vertexp{\Hy'}\in\Hy'$. In this case (see above),
every vertex of $\Hy'$ is a gamma leaf of $\Hy'$, and 
two edges are either comparable or disjoint. 
In particular, since $e_x\neq\Vertexp{\Hy}$ (by the previous point), 
there is a vertex $y$ such that any neighbour of $y$ is outside $e_x\smallsetminus\{x\}$,
therefore $y$ is also a gamma leaf of $\Hy$, not neighbour of $x$ in $\Hy\smallsetminus\{\Vertexp{\Hy}\}$.

In the other case, since $e_x\neq\Vertexp{\Hy}$ (by the previous point), 
we have $e_x\smallsetminus\{x\}\neq\Vertexp{\Hy[\smallsetminus\{x\}]}$.
The hypergraph $\Hy'$ is gamma acyclic and therefore, by induction, has two gamma leaves $y$ and $z$
that are not neighbours in $\Hy'\smallsetminus\Vertexp{\Hy'}$, and therefore not neighbours in $\Hy'$,
in particular in the edge $e_x\smallsetminus\{x\}$ in $\Hy'$. 
One of them, say $y$, does not belong to this latter edge, we call $e_y$ the maximal edge holding $y$ in $\Hy'$.
If $y$ is a singleton vertex, then it is a gamma leaf of $\Hy$, we assume it is not.
Let $e$ be the maximal edge holding $y$ in $\Hy'\smallsetminus\{e_y\}$.
Assume $e\cap e_x\neq \emptyset$.
Then $\Hy'$ is not gamma acyclic, a contradiction.
Finally, the neighbours of $y$ in $\Hy\smallsetminus\{e_y\}$ are beta leaves of $\Hy$,
hence $y$ is a gamma leave of $\Hy$.

We have proved the fact; we can apply it twice and get the expected result.
\end{proofen}
\end{thmen}
\subsection{Main Result}
\noindent
So far, we have proved:
\begin{description}
\item[1, easy to prove] If $\Hy'$ is obtained from $\Hy$ by removing a leaf, then $\Hy$ is acyclic iff $\Hy'$ is. (subsection 1.2)
\item[2, harder to prove] 
If $\Hy$ is acyclic but not empty, then it has a leaf. (subsection 2.1)
\end{description}
With this, the result can be obtained (see below).
But now, how will it help to prove other characterizations based on reduction processes?
We can keep the same pattern:
\begin{description}
\item[1, easy to prove] If $\Hy'$ is obtained from $\Hy$ by applying a step of the reduction process, then $\Hy$ is acyclic iff $\Hy'$ is. This is exactly done, for example, in \thref{lem:gammaeq}.
\item[2, harder to prove] If $\Hy$ is acyclic but not empty, 
	then one step of the reduction process can be applied.
\end{description}
In order to make the step 2 of the proof easier, we separate this into two parts:
\begin{description}
\item[2a, already done] If $\Hy$ is acyclic but not empty, then it has a leaf. (subsection 2.1)
\item[2b, easy to prove] If $\Hy$ has a leaf, 
	then one step of the reduction process can be applied. See, for example, \thref{lem:gammared}.
\end{description}
\begin{charen}[main result]
\label{char:hard}
Let $\Hy$ be a hypergraph.
The following two assertions are equivalent:
\begin{description}
\item[(1a)] The hypergraph $\Hy$ is alpha (resp. beta, gamma) acyclic.
\item[(2a)] The hypergraph $\Hy$ admits an alpha (resp. beta, gamma) elimination order.
\end{description}
\begin{proofen}
\thref{lem:easyway} states the (2a)$\Rightarrow$(1a) part.

We prove by induction on the number of vertices that an alpha (resp. beta, gamma) acyclic hypergraph
has an alpha (resp. beta, gamma) elimination order.
If a hypergraph has no vertex, then it is the empty hypergraph and the result is trivial.
Let us assume the fact holds for every hypergraph of less than $n$ vertices,
and consider $\Hy$ a hypergraph of $n$ vertices.

Assume $\Hy$ is alpha acyclic.
If $\Vertexp{\Hy}\notin\Hy$, \thref{thm:alphaleaves} asserts that it has alpha leaves.
If $\Vertexp{\Hy}\in\Hy$, then every vertex is an alpha leaf, by \thref{lem:trivialalpha}.

Assume $\Hy$ is beta (resp. gamma) acyclic.
If it has a single vertex, then this vertex is a beta (resp. gamma) leaf.
In the other case, \thref{thm:betaleaves} (resp. \thref{thm:gammaleaves}) proves $\Hy$ has beta (resp. gamma) leaves.

We have proved that if $\Hy$ is alpha (resp. beta, gamma) acyclic, then it has an alpha (resp. beta, gamma) leaf $x$.
By \thref{lem:eqleaf}, $\Hy[\smallsetminus\{x\}]$ is also alpha (resp. beta, gamma) acyclic, therefore,
by induction, $\Hy[\smallsetminus\{x\}]$ has an alpha (resp. beta, gamma) elimination order, therefore $\Hy$ also has one.
\end{proofen}
\end{charen}
We have proved that alpha (resp. beta, gamma) acyclicity is equivalent to the existence
of an alpha (resp. beta, gamma) elimination order; nevertheless, we have not
proved that removing repeatedly alpha (resp. beta, gamma) leaves is a confluent process,
i.e. given a starting hypergraph $\Hy$, the hypergraph $\Hy'$ obtained by removing
leaves until there is none left do not depend on the choices made: at each step, which leave to remove first?
Even if this process \emph{is actually confluent}, we do not need to prove so much (even if it would not be hard),
we only want to prove something weaker: reaching the empty hypergraph or not
does not depend on the choices made. This has the advantage of being
a straightforward consequence of the previous results.

Furthermore, we will be able to prove this weakened confluence property
for all other reduction processes. In order to avoid too much repetition,
we will not even state this property in the case of the reduction processes introduced section~3.
\begin{coren}[weakened confluence property]
\label{cor:confluent}
Let $\Hy$ be a hypergraph, having an alpha (resp. beta, gamma) leaf we call $x$.
Then $\Hy$ has an alpha (resp. beta, gamma) elimination order iff $\Hy[\smallsetminus\{x\}]$ has.
\begin{proofen}
By \thref{char:hard} ($\Hy$ is acyclic iff $\Hy$ has an elimination order) and \thref{lem:eqleaf} (if $x$ is a leaf of $\Hy$, then $\Hy$ is acyclic iff $\Hy[\smallsetminus\{x\}]$ is).
\end{proofen}
\end{coren}
\begin{coren}
Alpha, beta and gamma acyclicity are polynomial-time decidable.
\begin{proofen}
It is easy to see that checking
whether a hypergraph has an alpha (resp. beta, gamma) leaf, 
is polynomial; therefore,
by \thref{cor:confluent} (weakened confluence property), checking whether
a hypergraph has an alpha (resp. beta, gamma) elimination order is polynomial.
By \thref{char:hard}, this test decides alpha (resp. beta, gamma) acyclicity.
\end{proofen}
\end{coren}
Remark that, by \cite{TarjanY84}, alpha acyclicity is even linear-time decidable.

We state another corollary of these results, that illustrates the interest of the strengthening in:
``a acyclic hypergraph $\Hy$ has two leaves \emph{that are not neighbours in $\Hy\smallsetminus\{\Vertexp{\Hy}\}$}.
This result is in the same spirit as the \emph{sacred node property} in \cite{desi}.
\begin{thmen}[sacred node principle]
\label{thm:sacred}
Let $\Hy$ be an alpha (resp. beta, gamma) acyclic hypergraph, and let $e\in\Hy$.
Then $\Hy$ has an alpha (resp. beta, gamma) elimination order that eliminates every vertex ``outside'' $e$, i.e. in $\Vertexp{\Hy}\smallsetminus e$,
before removing any vertex of $e$.
\begin{proofen}[sketch]
Easy by induction. If $e$ is not a full edge of $\Hy$, then $\Hy$ has two leaves that are not neighbours in $\Hy\smallsetminus\{\Vertexp{\Hy}\}$,
in particular, one of them is not in $e$.
\end{proofen}
\end{thmen}
\begin{remen}
For the sake of completeness, we could introduce a notion of ``pure leaf'' and a related
notion of ``pure elimination order'' such that a hypergraph is cycle-free if and only if it has 
a pure elimination order.
It can be done as follows:
Let $\Hy$ be a hypergraph, and $x\in\Vertexp{\Hy}$ is \emph{a pure leaf} of $\Hy$ if:
\begin{itemize}
\item We cannot find $y$ and $z$ such that $\Mp{\Hy[\{x,y,z\}]}=\{\{x,y\},\{y,z\},\{x,z\}\}$.
\item The neighbourhood of $x$ in $\Hy$ is a clique of $\Hy$, i.e. for any $y$ and $z$ in $\Vertexp{\Hy(x)}$,
some edge $e\in\Hy$ satisfies $\{y,z\}\subseteq e$.
\end{itemize}
As we did section 1.2, it is easy, by using the facts of \thref{lem:trivialfacts},
 to prove that, if $x$ is a pure leaf of $\Hy$,
then $\Hy$ is cycle-free iff $\Hy[\smallsetminus\{x\}]$ is.
There remains to prove that a cycle-free hypergraph has a pure leaf.

So, let $\Hy$ be a cycle-free hypergraph.
Consider:
\mydisplaystyle{\begin{align*}
\Hy'=\Hy\cup\{K\subseteq\Vertexp{\Hy}\,|\,K\ \mathrm{is\ a\ maximal\ clique\ of}\ \Hy\}
\end{align*}}
Obviously, two vertices are neighbours in $\Hy$ iff they are neighbours in $\Hy'$.
As a consequence, $\Hy$ and $\Hy'$ have the same cliques.

Now observe that $\Hy'$ is alpha acyclic.
Therefore, by \thref{thm:alphaleaves}, it has an alpha leaf $x$.
Trivially, the neighbourhood of $x$ in $\Hy'$ is a clique of $\Hy'$, hence the neighbourhood of $x$ in $\Hy$ it is a clique of $\Hy$,
 and, since $\Hy$ is cycle-free, $x$ is not in a triangle.
We have proved $x$ is a pure leaf of $\Hy$.
\end{remen}
\section{Other Characterizations}
\noindent
Here we make the connection with other known characterizations of alpha and gamma acyclicity,
and show how our notions of leaf allow to prove easily their equivalence.
The proofs of this section are not as detailed as the other proofs of this document,
because the notions of join tree for example are in fact rather heavy.

We define the operations needed to describe the different considered reduction processes.
\begin{defen}[operations]
We define the following operations:
\begin{description}
\item[linearization]
	We say $\Hy'$ is obtained by $\Hy$ by \emph{linearization} iff
	$\Hy'=\Hy[\smallsetminus\{x\}]$ where $\exists y\in\Vertexp{\Hy},y\neq x\:\Hy(x)=\Hy(y)$, i.e. $x$ and $y$ are contained in exactly the same edges
	of $\Hy$.
\item[singleton edge removal]
	We say $\Hy'$ is obtained by $\Hy$ by \emph{singleton edge removal} iff
	$\Hy'=\Hy\smallsetminus\{e\}$ where $e\in\Hy$ and $\mathrm{card}\,e=1$, i.e. $e$ is in the form $\{x\}$.
\item[included edge removal]
	We say $\Hy'$ is obtained by $\Hy$ by \emph{included edge removal} iff
	$\Hy'=\Hy\smallsetminus\{e\}$ where $e\in\Hy$ and $\exists f\in\Hy\:e\subset f$.
\item[singleton vertex removal] By analogy with singleton edge removal, 
	we say $\Hy'$ is obtained by $\Hy$ by \emph{singleton vertex removal} iff
	$\Hy'=\Hy[\smallsetminus\{x\}]$ where $\mathrm{card}\,\Hy(x)=1$, i.e. $x$ was contained in a single edge of $\Hy$.
\end{description}
\end{defen}
\begin{remen}
For some of these operations, applying it until not possible defines a hypergraph transformation
that we have already introduced. For example, 
applying ``included edge removal'' of a hypergraph $\Hy$ until not possible leads to the hypergraph $\M(\Hy)$.
Applying ``linearization'' of a hypergraph $\Hy$ until not possible, also known as ``contracting all modules'' in the literature,
leads to a hypergraph that is isomorphic to the normalization of $\Hy$ denoted $\N(\Hy)$ in this document.

As a consequence, properties that are invariant w.r.t. $\M$ (resp. $\N$) are also invariant w.r.t. included edge removal (resp. linearization).
Let us prove it briefly.

Let $P$ a property such that, for all $\Hy$, we have $P(\Hy)\Leftrightarrow P(\M(\Hy))$.
Let $\Hy'$ obtained from $\Hy$ by included edge removal.
We have $\M(\Hy')=\M(\Hy)$, therefore $P(\Hy)\Leftrightarrow P(\M(\Hy))\Leftrightarrow P(\M(\Hy'))\Leftrightarrow P(\Hy')$.

As a corollary, by \thref{lem:stabinduc}, 
cycle-freedom, conformity and alpha acyclicity are invariant w.r.t. included edge removal.
\end{remen}
\subsection{Other Alpha Acyclicity Characterizations}
\subsubsection{GYO Reduction}
\noindent
In this section, we make the connection between alpha elimination order and GYO reduction.
This is the original GYO (Graham, Yu, \"Ozsoyoglu)
non deterministic algorithm described by 
\cite{Graham,YO}
(see also \cite{AbiteboulHV95} ex. 5.29 p151 and 
\cite{DBLP:journals/tods/FaginMU82}).
\begin{defen}
The \emph{Graham-Yu-\"Ozsoyoglu} (\emph{GYO}, for short) \emph{operations} are the following:
\begin{itemize}
\item included edge removal, and
\item singleton vertex removal.
\end{itemize}
A hypergraph $\Hy$ is \emph{GYO-reducible} if there exists
a sequence of GYO operations that leads to the empty hypergraph.
\end{defen}
\begin{remen}
Recall that our minimization operation $\M$ is the fixed-point closure of ``included edge removal''.
\end{remen}
\begin{charen}
A hypergraph $\Hy$ is alpha acyclic iff:
\begin{description}
\item[($\alpha$2b)] $\Hy$ is GYO-reducible.
\end{description}
\begin{proofen}
By induction on the hypergraph \emph{size}. The equivalence is trivial for the empty hypergraph.
Assume this holds for every hypergraph of size less than $n$.
Take $\Hy$ of size $n$.

Assume $\Hy$ is alpha acyclic, therefore it has an elimination order, therefore it has a leaf $x$.
If $x$ is a singleton vertex, then we can apply singleton vertex removal to $\Hy$, and 
by induction $\Hy[\smallsetminus\{x\}]$ is GYO-reducible.
Assume $x$ is not a singleton vertex.
As an alpha leaf, $x$ is a singleton vertex in $\M(\Hy)$, therefore $\Hy\neq\M(\Hy)$, therefore we can apply
included edge removal and get $\Hy'$. Since $\M(\Hy')=\M(\Hy)$, we can deduce, by \thref{lem:stabinduc},
that $\Hy'$ is acyclic,
and GYO-reducible by induction.

Now we assume $\Hy$ is GYO-reducible.
If the first operation is an included edge removal, then the hypergraph obtained after this step is acyclic by the induction hypothesis,
then, by \thref{lem:stabinduc} again, so is $\Hy$.
If the first operation is a singleton vertex $x$ removal, then $x$
is a singleton vertex in $\Hy$, and by induction $\Hy[\smallsetminus\{x\}]$ is acyclic, therefore
has an alpha elimination order.
Since $x$ is a leaf, $\Hy$ has an elimination order and is therefore acyclic.
\end{proofen}
\end{charen}
\subsubsection{Join Trees}
\noindent
In this section, we make the connection between GYO reduction and join tree (also called \emph{junction tree}), 
a very common characterization.
See for example \cite{Robertson1986309}, or the seminal paper \cite{yann} for databases, or \cite{rina} for constraint satisfaction problems.
A join tree can be seen as a condensed representation of a \emph{realized} GYO elimination process.
\begin{defen}
A \emph{tree} is an acyclic graph. A \emph{join tree} is a couple $(\mathcal{T},L)$ where
$\mathcal{T}$ is a tree and $L$ is a labeling function defined on the domain $\Vertexp{\mathcal{T}}$,
that satisfies the so-called \emph{join property}. The \emph{join property} states that,
for all vertices $x$ and $y$ such that $L(x)\cap L(y)\neq \emptyset$, there is a path
from $x$ to $y$, and each vertex $v$ on the path satisfies $L(x)\cap L(y)\subseteq L(v)$.

A join tree $(\mathcal{T},L)$ is \emph{the join tree of a hypergraph $\Hy$} when
$L$ is injective and the set of images of $L$ is exactly $\Hy$.
A hypergraph $\Hy$ \emph{has a join tree} 
when there exists a join tree of $\Hy$.
\end{defen}
The following theorem is very classical.
We only make a short proof, a full proof of it can be found in \cite{Acharya} for example.
\begin{charen}[join tree]
\label{thm:jt}
The hypergraph $\Hy$ is alpha acyclic iff:
\begin{description}
\item[(2c)] The hypergraph $\Hy$ has a join tree.
\end{description}
\begin{proofen}[sketch]
Consider the following fact:
\begin{description}
\item[fact] If $\Hy'$ is obtained from $\Hy$ by a GYO operation,
then $\Hy'$ has a join tree iff $\Hy$ has a join tree.
\end{description}
Using this fact inside a proof by induction, 
together with the equivalence between GYO-reducible and alpha acyclic,
allows to establish the expected result.

Now, we prove the non-trivial part of the fact, i.e. 
if $\Hy'$ is obtained from $\Hy$ by a GYO operation and $\Hy'$ has a join tree, then $\Hy$ also has one.

If $\Hy'$ is obtained by included edge removal, 
then the edge $e\in\Hy\smallsetminus\Hy'$ is included in some other edge $f\in\Hy'$.
Let $x$ be the vertex in $\mathcal{T}$ such that $L(x)=f$.
Take a fresh symbol $y\notin\Vertexp{\mathcal{T}}$, and
extend the definition of $L$ with $L':y\mapsto e$. 
It is easy to check that $(\mathcal{T}\cup\{\{x,y\}\},L')$
satisfies the join property, and is therefore a join tree of $\Hy$.

In the other case, $\Hy'$ is obtained by removing a singleton vertex we call $t$.
There is only one edge $e_t$ that includes it. 
Notice $e_t\smallsetminus\{t\}\in\Hy'$,
let $x$ be the vertex of $\mathcal{T}$ such that $L(x)=e_t\smallsetminus\{t\}$.
If $e_t\smallsetminus\{t\}$ does not belongs to $\Hy$, the result is obvious.

In the other case, take $y\notin\Vertexp{\mathcal{T}}$,
extend $L$ into $L'$ by adding $L'(y)=e_t$.
It is easy to check that $(\mathcal{T}\cup\{\{x,y\}\},L')$
is a join tree of $\Hy$.
\end{proofen}
\end{charen}
\begin{charennp}[alpha acyclicity]
\label{char:alpha}
Let $\Hy$ be a hypergraph.
The following are equivalent:
\begin{description}
\item[($\alpha$1a)] The hypergraph $\Hy$ is alpha acyclic, i.e. conformal and cycle-free.
\item[($\alpha$1b)] We \emph{cannot} find $S\subseteq\Vertexp{\Hy}$ such that $\M(\Hy[S])$ is either $\{S\smallsetminus\{x\}\,|\,x\in S\}$ or a usual graph cycle.
\item[($\alpha$2a)] The hypergraph $\Hy$ admits an alpha elimination order.
\item[($\alpha$2b)] The hypergraph $\Hy$ is GYO-reducible.
\item[($\alpha$2c)] The hypergraph $\Hy$ has a join tree.
\end{description}
\end{charennp}
\subsection{Beta Acyclicity}
\noindent
Notice that beta acyclicity does not seem to admit a simple characterization in terms
of join tree, by contrast with alpha and gamma acyclicity.
\begin{charennp}[beta acyclicity]
\label{char:beta}
Let $\Hy$ be a hypergraph.
The following are equivalent:
\begin{description}
\item[($\beta$1a)] The hypergraph $\Hy$ is beta acyclic, i.e. every subset of $\Hy$ is alpha acyclic.
\item[($\beta$1b)] Every subset of $\Hy$ is cycle-free.
\item[($\beta$1b)] (reformulated) We cannot get a usual graph cycle from $\Hy$ by removing vertices and/or edges.
\item[($\beta$2a)] The hypergraph $\Hy$ admits a beta elimination order.
\item[($\beta$2b)] (reformulation of ($\beta$2a)) The hypergraph $\Hy$ can be reduced to the empty set by repeatedly removing nest points.
\end{description}
\end{charennp}
These characterizations are ``all we have'' for beta acyclicity;
nonetheless, it is not a problem if we consider how simple and natural the characterizations ($\beta$1b) and ($\beta$2a) are.
\begin{remen}
Characterization ($\beta$2a) admits as a corollary\footnote{More precisely: this characterization gives immediately that $u_n=u_{n-1}+n$ where $u_n$
is the maximal number of edges of a beta acyclic hypergraph with $n$ vertices. We can easily prove the result by induction.} 
that a beta acyclic hypergraph on $n$ vertices has at most \smash{$\frac{n(n+1)}{2}$} edges.
Observe that the maximum number of edges is reached in the case of an interval hypergraph, i.e. $\{\{x_i,\ldots,x_j\}\,|\,i\leq j\leq n\}$.

It is obvious that gamma acyclicity share this bound (it is a restriction of a beta acyclicity).
By contrast, alpha acyclicity have no bound: 
for any $n$, 
the ``full hypergraph'', that is: \smash{$\big\{e\subseteq\{1,\ldots,n\}\,\big|\,e\neq\emptyset\big\}$}, is alpha acyclic.
\end{remen}
\subsection{Other Gamma Acyclicity Characterizations}
\begin{defen}[DM reduction]
The \emph{D'Atri-Moscarini} (\emph{DM}, for short) \emph{operations} are the following:
\begin{itemize}
\item singleton vertex removal,
\item singleton edge removal, and
\item linearization.
\end{itemize}
A hypergraph is \emph{DM-reducible} iff there exists
a sequence of DM operations that leads to the empty hypergraph.
\end{defen}
First of all, we state that gamma acyclicity is invariant w.r.t. DM operations.
\begin{lemen}
\label{lem:gammaeq}
Let $\Hy$ be a hypergraph. For any hypergraph $\Hy'$ obtained from $\Hy$ by applying one DM operation, $\Hy$ is gamma acyclic if and only if $\Hy'$ is.
\begin{proofen}
Notice that \thref{lem:eqleaf} already states the result in the case of \emph{singleton vertex removal}.
Proving the result in the cases of \emph{singleton edge removal} and of \emph{linearization} is easy and left to the reader;
using the characterization ($\gamma$1b) of gamma acyclicity (see \thref{char:type1}) makes it easier.
\end{proofen}
\end{lemen}
Now we give a simple lemma that allows to make the connection
between a gamma leaf and the reduction process from \cite{82-DaMo} that also characterizes gamma acyclicity.
As mentioned just before stating the main result (section 2.2), we prove that, when there is a gamma leaf, then some operation of the gamma reduction process,
the DM-reduction, can be applied.
\begin{lemen}
\label{lem:gammared}
If a hypergraph has a gamma leaf, then at least one step of the DM reduction process can be performed.
\begin{proofen}
We prove the following fact.
A gamma leaf of a hypergraph $\Hy$ either:
\begin{enumerate}
\item has a neighbour that is a singleton vertex, 
\item has a neighbour in a singleton edge, or
\item has two distinct neighbours $x$ and $y$ such that $\Hy(x)=\Hy(y)$.
\end{enumerate}
Let be $t$ a gamma leaf of $\Hy$, assume (1) and (2) are false.
We call $e_t$ the maximal edge containing $t$.
By definition of a gamma leaf, any neighbour of $t$ in $\Hy\smallsetminus\{e_t\}$
is also a gamma leaf.
Let $x$ be, among them, a vertex contained in a maximal number of edges.

Consider the smallest edge $e_1$ and the biggest edge $e_2=e_t$ containing $x$.
Since hypothesis (1) is false, the edge $e_1$ cannot hold only $x$,
therefore $e_1$ contains at least another vertex we call $y$. 
Since hypothesis (2) is false, $e_1\neq e_2$, therefore $x$ and $y$ are neighbours in $\Hy\smallsetminus\{e_2\}$.
By definition of a gamma leaf,
every neighbour of $x$ in $\Hy\smallsetminus\{e_2\}$ is a beta leaf of $\Hy$, therefore $y$
is also a beta leaf of $\Hy$.

Assume some edge $e$ includes $x$ but not $y$. 
Since $x$ is a beta leaf, we must have $e\subseteq e_1$, which contradicts the definition of $e_1$.
Assume some edge $e$ includes $y$ but not $x$.
Since $y$ is a beta leaf, we have $e\subseteq e_1$.
Then $y$ belongs to every edge that includes $e_1$ and in $e$, and also belongs to $e$
therefore it belongs to more edges than $x$ does, a contradiction.

We have proved $\Hy(x)=\Hy(y)$, therefore the fact is proved.
Now we can conclude. Assume $t$ be a gamma leaf of $\Hy$.
If $t$ has a neighbour that is a singleton vertex, then we can apply \emph{singleton vertex removal}.
If $t$ has a neighbour in a singleton edge, then we can apply \emph{singleton edge removal}.
If $t$ has two distinct neighbours $x$ and $y$ such that $\Hy(x)=\Hy(y)$, we can apply \emph{linearization}.
This concludes the proof.
\end{proofen}
\end{lemen}
\begin{charen}
A hypergraph $\Hy$ is gamma acyclic iff:
\label{charDM}
\begin{description}
\item[($\gamma$2b)] The hypergraph $\Hy$ is DM-reducible.
\end{description}
\begin{proofen}
We prove a hypergraph $\Hy$ is DM-reducible iff it is gamma-acyclic, by induction on the size of $\Hy$.
Let $\Hy$ be a hypergraph of size $n$.

If $\Hy$ is DM-reducible, we call $\Hy'$ the hypergraph obtained after one step of the DM-reduction.
This hypergraph $\Hy'$ is also DM-reducible hence, by induction, gamma acyclic.
By \thref{lem:gammaeq}, so is $\Hy$.

If $\Hy$ is gamma acyclic, then it has a gamma leaf.
By~\thref{lem:gammared}, one of the operations of the DM-reduction
--- i.e. linearization, singleton edge removal or singleton vertex removal ---
can be performed, we call $\Hy'$ the
resulting hypergraph, which is gamma acyclic by \thref{lem:gammaeq}.
By induction, $\Hy'$ is DM-reducible, therefore so is $\Hy$.
\end{proofen}
\end{charen}
%
\begin{charennp}[gamma acyclicity]
\label{char:gamma}
Let $\Hy$ be a hypergraph.
The following are equivalent:
\begin{description}
\item[($\gamma$1a)] The hypergraph $\Hy$ is gamma acyclic, i.e. $\Hy$ is beta acyclic and we cannot find $x,y,z$ such that $\{\{x,y\},\{x,z\},\{x,y,z\}\}\subseteq \Hy[\{x,y,z\}]$. (our definition)
\item[($\gamma$1b)] The hypergraph $\Hy$ is cycle-free and we cannot find $x,y,z$ such that $\{\{x,y\},\{x,z\},\{x,y,z\}\}\subseteq \Hy[\{x,y,z\}]$. 
(Definition~3 of gamma acyclicity in \cite{Fagin83degreesof})
\item[($\gamma$2a)] The hypergraph $\Hy$ admits a gamma elimination order.
\item[($\gamma$2b)] The hypergraph $\Hy$ is DM-reducible.
\item[($\gamma$2c)] \cite{theseDuris,betagamma} For any $e\in\Hy$, $\Hy$ has a rooted join tree with disjoint branches\footnote{%
A rooted join tree with disjoint branches, in addition to the join property,  
satisfies: for any vertices of the tree $a$ and $b$ with no ancestor relationship between them, the label of $a$ and the label of $b$ are disjoint sets.} 
whose root is labelled $e$. The equivalence between this characterization and ($\gamma$2b) is rather easy to prove\rlap{\footnotemark}.
\footnotetext{More precisely:
We can very easily adapt \thref{lem:gammared} to prove that, if $\Hy$ is a gamma acyclic hypergraph,
then we may perform two steps of the DM reduction process on non-neighbour vertices 
by the ``sacred node principle'', i.e. \thref{thm:sacred}, therefore we can choose a sequence of DM operations that will
preserve a given edge $e$ until there is only $e$ left in $\Hy$.
It will be easy to show how to build a join tree rooted in $e$ with disjoint branches with this DM reduction.}
\end{description}
\end{charennp}
\section{Relevancy of Acyclicity Notions}
\noindent
We aim at introducing a notion of ``good acyclicity notion'',
in order to check whether new interesting acyclicity notions
could be introduced or not.

Note: this section makes heavy use of ``closed under'' and ``invariant w.r.t.''; 
their respective definitions can be found in \thref{def:prop}, page~\pageref{def:prop}.
\subsection{Reasonable Hypergraph Acyclicity Notions}
\noindent
First of all, we start with a counter-example:
so far, we did not mention Berge acyclicity.
In Fagin's paper, this notion of Berge acyclicity was also discarded,
with the argument that the notion was too restrictive.
We give a totally different argument.

Berge acyclicity admits a very natural definition \cite{Berge:1985}:
\begin{itemize}
\item A hypergraph $\Hy$ is Berge acyclic when the graph $\mathcal{G}=\big\{\{x,e\}\,\big|\,x\in e\,\mathrm{and}\,e\in\Hy\big\}$ is acyclic. 
\end{itemize}
For example, the hypergraph 
$\Hy=\{\{x,y\},\{x,y,z\}\}$, 
that we call ``Berge triangle'' in Fig.~\ref{fig:exhyper},
is not Berge acyclic: 
\mydisplaystyle{\begin{align*}
\mathcal{G}=\{\{x,e\},\{y,e\},\{x,f\},\{y,f\},\{z,f\}\}
\end{align*}}
 is not gamma acyclic.
Notice the edge $\{z,f\}$ is not responsible for the graph $\mathcal{G}$ not being acyclic.
Is there a simpler hypergraph that is not acyclic?
The answer is no: the vertex $z$ is needed to make a difference between the two edges of the hypergraph.

But now, if we consider \emph{multi}-hypergraphs, i.e. hypergraphs where there may be several edges that contain the same set of vertices,
then the multi-hypergraph $\mathcal{M}=[\{x,y\},\{x,y\}]$ is not Berge acyclic, because
$\mathcal{G}=\{\{x,e\},\{y,e\},\{x,f\},\{y,f\}\}$ is not acyclic.
This suggests that Berge acyclicity is not an actual \emph{hypergraph notion}, 
but rather a \emph{multi-hypergraph notion}; see the remark after the definition.

Now we introduce a more formal criterion that a ``reasonable hypergraph acyclicity property'' should satisfy.
Intuitively, we would like a multi-hypergraph notion that ignores the number of copies the ``same edge''.
Dually, we also require that it ignores the numbers of copies of the ``same vertex''.
This leads to require a hypergraph property $P$ (that ignores, by definition of a hypergraph, the number of copies of a given edge)
such that $P(\Hy)\Leftrightarrow P(\mathcal{N}(\Hy))$, i.e. that ignores the number of ``copies of'' a vertex.
\begin{defen}
A property is a \emph{reasonable hypergraph property} when it is invariant w.r.t. normalization $\mathcal{N}$.
That is to say: if $\Hy$ is a hypergraph such that $\Hy(x)=\Hy(y)$, then a reasonable hypergraph property $P$
should satisfy $P(\Hy)\Leftrightarrow P(\Hy[\smallsetminus\{y\}])$.

A property $P$ is a \emph{reasonable hypergraph acyclicity notion} (or, for short:\emph{reasonable acyclicity notion})
when:
\begin{enumerate}
\item It is a reasonable hypergraph property.
\item For every \emph{graph} $\mathcal{G}$, $P(\mathcal{G})$ if and only if $\mathcal{G}$ is a acyclic graph.
\item The property $P$ is invariant w.r.t. singleton edge addition and deletion.
\item The property $P$ is invariant w.r.t. singleton vertex addition and deletion.
\end{enumerate}
These properties will be referred to as (1), (2), etc.
\end{defen}
\begin{remen}
Instead of requiring that the property $P$ coincides with graph acyclicity on all graphs, 
we could \emph{only} require that $P$ coincides with the usual graph acyclicity notion on non-acyclic graphs and on the empty hypergraph. 
By using this property and the properties (3) and (4), 
we could then deduce that $P$ must also coincide with graph acyclicity on acyclic graphs.
\end{remen}
We will not always use all these properties; most results only need a few of them.
But the conjunction of all these properties,
which is not very constraining,
allows to consider a single definition of ``reasonable hypergraph acyclicity notion''.

Berge acyclicity is not a reasonable hypergraph acyclicity notion: it is not
a reasonable hypergraph property, and the condition (4) is not satisfied;
the ``Berge triangle'' gives a counter example in both cases.
Now if we consider multi-hypergraphs, the condition (4) is satisfied.
\begin{thmen}
Alpha, beta, and gamma acyclicity and cycle-freedom are reasonable hypergraph acyclicity notions.
\begin{proofen}[sketch]
Each of these proof is simple, but is always a not-so-short case study.
Notice that \thref{lem:eqleaf} already establishes that the notions have property (4).
\end{proofen}
\end{thmen}
We defined a reasonable hypergraph property as something that is invariant w.r.t. normalization.
Intuitively, in the dual of a normalized hypergraph, two different vertices will give rise to two different edges;
and taking the dual twice have no effect.

Furthermore, the criteria (3) and (4) are dual to each other.
This raises the question whether there is some reasonable hypergraph acyclicity property $P$ that is self-dual, i.e. $P(\Hy)\Leftrightarrow P(\mathcal{D}(\Hy))$.
The following theorem answers this question.
\begin{thmen}[self-dual acyclicity notions]
Beta and gamma acyclicity are invariant under duality, but alpha acyclicity and cycle-freedom are not. \label{thm:self-dual}
\begin{proofen}[sketch]
We prove the following: if $\mathcal{D}(\Hy)$ is beta (resp. gamma) acyclic, then so is $\Hy$.
Assume $\Hy$ is not beta acyclic. Then we can obtain a usual graph cycle from $\Hy$ by edges and vertices removal.
We therefore can obtain the dual of this usual graph cycle from $\mathcal{D}(\Hy)$ by vertices and edge removal.
This latter hypergraph is also a usual graph cycle.
Now we assume $\Hy$ is not gamma acyclic. If it is not beta acyclic, then we have the result by previous point.
In the other case, we can obtain a hypergraph isomorphic to the ``gamma triangle'', i.e. $\{\{x,y\},\{y,z\},\{x,y,z\}\}$ by edges and vertices
removal. We therefore can obtain the dual of this ``gamma triangle'' $\mathcal{D}(\Hy)$ by vertices and edge removal.
The latter is isomorphic to its dual.

We know that if $\mathcal{D}(\Hy)$ is beta (resp. gamma) acyclic, then so is $\Hy$.
Since beta and gamma acyclicity are invariant w.r.t. normalization (remind that $\mathcal{N}(\Hy)=\mathcal{D}(\mathcal{D}(\Hy))$), 
We have $\Hy$ is beta (resp. gamma acyclic) iff $\mathcal{D}(\mathcal{D}(\Hy))$ is, which implies $\mathcal{D}(\Hy)$ is beta (resp. gamma) acyclic
by the previous point. We have proved $\Hy$ is beta (resp. gamma) acyclic iff $\mathcal{D}(\Hy)$ is.

The hypergraph 
$\{\{x,y\},\{y,z\},\{x,z\},\{x,y,z\}\}$
is alpha acyclic (and therefore cycle-free) but its dual: 
$\{\{e,g,h\},\{e,f,h\},\{f,g,h\}\}$
is not cycle-free (and therefore not alpha acyclic).
\end{proofen}
Notice that even graph acyclicity does not have this interesting property of being closed under duality:
the dual of $\{\{x,t\},\{y,t\},\{z,t\}\}$ is not a graph.
\end{thmen}
\subsection{Desirable Closure Properties of Acyclicity Notions}
\noindent
In previous subsection,
we have introduced a definition of a ``reasonable hypergraph acyclicity notion'',
that was made as little constraining as possible.
Now, we will define a ``\emph{good} hypergraph acyclicity notion''
as a reasonable hypergraph acyclicity notion that enjoys
nice closure properties.

\noindent
%
A very well-known notion on graphs is the notion of \emph{graph minor} (see \cite{Diestel,Berge69,BM08,CLZ,Lov06});
many interesting graph properties are closed under minor, e.g. planarity.
Now we introduce the operations that define the notion of hypergraph minor reported in \cite{Duchet},
that is a well quasi-ordering by a result of Robertson and Seymour in 1987, published much later \cite{RobertsonS10}.
Intuitively, a ``good acyclicity notion'' should be closed under taking minor, i.e. it should
be closed under each of the following operations.
\begin{defen}
We define (or recall):
\begin{description}
\item[edge shrinking] Replace an edge by a subset of it.
\item[vertex removal] Remove a vertex.
\item[edge removal] Remove an edge.
\item[edge contraction] Take two neighbours $x$ and $y$, replace in the edges every occurrence of $y$ by $x$.
	That is to say a hypergraph $\Hy'$ is obtained by \emph{edge contraction} from $\Hy$ when there are two vertices $x$ and $y$ neighbours
	in $\Hy$ such that:
\mydisplaystyle{\begin{align*}
\Hy'=\{e\smallsetminus\{x\}\cup\{y\}\,|\,e\in\Hy(x)\}\cup(\Hy\smallsetminus\Hy(x))
\end{align*}}
\end{description}
We say a hypergraph $\Hy$ is a \emph{minor} of another hypergraph $\Hy'$ if $\Hy$ is isomorphic to a hypergraph obtained from $\Hy'$
by applying one or several of the operations listed above, i.e. edge shrinking, vertex removal, edge removal, and edge contraction.
\end{defen}
\begin{remen}
Observe that Berge acyclicity is closed under each of these operations.
This proves informally that it is a ``good multi-hypergraph acyclicity notion''.
\end{remen}
We would like to prove that alpha, beta and gamma acyclicity are ``good acyclicity notions''
with the following argument: they are reasonable hypergraph acyclicity notions that are 
closed under taking hypergraph minor.
Unfortunately, this is not the case. Even worse, no such acyclicity notion can even exist: 
\begin{thmen}[limits of closure properties]
\label{thm:limclosure}
A reasonable hypergraph acyclicity property cannot be:
\begin{itemize}
\item closed under edge shrinking.
\item closed under both edge removal and edge contraction.
\end{itemize}
\begin{proofen}
Let $P$ be a reasonable hypergraph acyclicity.

Consider \smash{$\Hy_1=\{\{x,a,b,c\},\{y,a,b,c\},\{z,a,b,c\}\}$}.
Let \smash{$\Hy_2=\mathcal{N}(\Hy_1)$}, i.e. \smash{$\Hy_2=\{\{x,a\},\{y,a\},\{z,a\}\}$}.
By property (2), \smash{$P(\Hy_2)$} is true; by property (1), therefore so is \smash{$P(\Hy_1)$}.
Now let \smash{$\Hy_3=\{\{a,b\},\{b,c\},\{a,c\}\}$}.
By property (2), \smash{$P(\Hy_3)$} is false.
We can notice \smash{$\Hy_3$} is obtained from \smash{$\Hy_1$} by edge shrinking,
but we have proved \smash{$P(\Hy_1)$} is true and \smash{$P(\Hy_3)$} is not, therefore $P$ is not closed under
edge shrinking.

Consider \smash{$\Hy_1=\{\{x\}\}$}. By (2), \smash{$\Hy_1$} satisfies $P$, and therefore,
by (1), so does \smash{$\Hy_2=\{\{x,y,z\}\}$}. 
By (3), so does \smash{$\Hy_3=\{\{x,y,z\},\{x\},\{y\},\{z\}\}$}.
By (1), so does 
\smash{$\Hy_4=\{\{x_1,x_2,y_1,y_2,z_1,z_2\},\{x_1,x_2\},\{y_1,y_2\},\{z_1,z_2\}\}$}.

We can obtain 
\smash{$\Hy_5=\{\{x_1,y_1,z_1\},\{x_1,y_1\},\{y_1,z_1\},\{z_1,x_1\}\}$} from \smash{$\Hy_4$}
by edge contraction: choose \smash{$x_2$} and \smash{$y_1$}, then \smash{$y_2$} and \smash{$z_1$}, and, finally, \smash{$z_2$} and \smash{$x_1$}.
We can obtain \smash{$\Hy_6=\{\{x_1,y_1\},\{y_1,z_1\},\{z_1,x_1\}\}$} from \smash{$\Hy_4$} by edge removal.
We have obtained a hypergraph \smash{$\Hy_6$}, that does not satisfy $P$ by (2),
by edge contraction and edge removal from \smash{$\Hy_4$} that does satisfy $P$.
Therefore $P$ is not closed under both edge contraction and edge removal.
\end{proofen}
\end{thmen}
No reasonable hypergraph acyclicity notion can be closed under taking minor!
Instead of defining a good acyclicity notion as a reasonable hypergraph acyclicity notion
that is closed under taking minor, we define it in a relaxed setting:
\begin{defen}[good acyclicity notion]
A \emph{good acyclicity notion} is a reasonable hypergraph acyclicity notion that is closed under
two operations among the three operations \emph{vertex removal}, \emph{edge removal}, and \emph{edge contraction}.
\end{defen}
Notice that, since \thref{thm:limclosure} states that no reasonable acyclicity notion can be
closed under both edge removal and edge contraction, a good acyclicity notion is necessarily closed under
vertex removal.
\begin{remen}
In \cite{DBLP:journals/tcs/AdlerGK12}, they introduce an \emph{ad hoc} notion of hypergraph minor,
such that alpha acyclicity is closed under this new definition,
but neither beta nor gamma acyclicity are closed under taking minor with this definition.
\end{remen}
\thref{lem:stabinduc} already states that cycle-freedom, alpha acyclicity, beta acyclicity, 
and gamma acyclicity are closed under vertex removal, 
and that beta and gamma acyclicity are closed under edge removal.
Let us prove additional closure properties.

One of the following results is not new: the paper \cite{DBLP:journals/tcs/AdlerGK12} already establishes
that having hypertree depth $k$ is closed under edge contraction, 
which implies trivially that alpha-acyclicity is closed under edge contraction.
\begin{thmen}[properties closed under edge contraction]
\label{thm:contraction}
Cycle-freedom and alpha acyclicity are closed under edge contraction.
\begin{proofen}[easy but long]
Let $\Hy$ be a cycle-free hypergraph.
Consider $\Hy'$ obtained by edge contraction of $x$ and $y$ in the edge $e$ of $\Hy$, i.e.
$\Hy'=\{e\smallsetminus\{x\}\cup\{y\}\,|\,e\in \Hy(x)\}\cup(\Hy\smallsetminus\Hy(x))$.
Assume $\Hy'$ is not cycle-free. Let $S$ such that $\Hy'[S]$ is a usual graph cycle.
If $y\notin S$ , then $\Hy[S]$ is the same usual graph cycle, therefore $\Hy$ is not cycle-free, a contradiction.
From now, we assume $y\in S$. We have 
$(x_1\ldots,x_{i-1},x_i=y,x_{i+2},\ldots,x_k)$
a cycle of $\Hy'[S]$.

Consider $\Hy''=\Hy[S\cup\{x\}]$.
The vertex $y$ has exactly two neighbours $x_{i-1}$ and $x_{i+1}$ in $\Hy'[S]$.
If $x_{i-1}$ and $x_{i+1}$ are neighbours of $y$ (resp. $x$) in $\Hy''$, then $\Hy''[\smallsetminus\{x\}]$ (resp. $\Hy''[\smallsetminus\{y\}]$) is a cycle, a contradiction.
From now we assume $x_{i-1}$ is neighbour or $y$ and $x_{i+1}$ is neighbour of $x$ in $\Hy''$.
It is now easy to see that $(x_1\ldots,x_{i-1},x_i=y,x,x_{i+2},\ldots,x_k)$ is a cycle of $\Hy''$, which is therefore
no cycle-free, which is a contradiction.

We have proved $\Hy'$ is cycle-free, therefore cycle-freedom is closed under edge contraction.

Let $\Hy$ be an alpha acyclic hypergraph.
Consider $\Hy'$ obtained by edge contraction of $x$ and $y$ in the edge $e$ of $\Hy$, i.e.
$\Hy'=\{e\smallsetminus\{x\}\cup\{y\}\,|\,e\in \Hy(x)\}\cup(\Hy\smallsetminus\Hy(x))$.
By the previous point, $\Hy'$ is cycle-free. Assume it is not conformal.
Let $S$ such that $\M(\Hy'[S])=\{S\smallsetminus\{x\}\,|\,x\in S\}$.
If $y\notin S$, then $\Hy[S]$ is a non-conformal clique, contradiction. 
If $y\in S$, then, in $\Hy''=\Hy[S\cup\{x\}]$, every vertex $z$ is neighbour of $x$
or neighbour of $y$.
Assume every vertex is neighbour of both. Then $S\cup\{x\}$ is a clique. If $S\cup\{x\}\in\Hy''$,
then $S\in\Hy'[S]$, which is a contradiction.

Then, for any vertex $z$, $z$ is neighbour of either $x$ or $y$ but not both.
It is easy to see that at least one vertex $u$ (resp. $v$) other than $y$ (resp. $x$) is neighbour of $x$ (resp. $y$).
Then $\M(\Hy''[\{x,y,u,v\}])=\{\{x,y\},\{y,v\},\{v,u\},\{u,x\}\}$ which is not cycle-free, a contradiction.
\end{proofen}
\end{thmen}
We have proved that alpha, beta, and gamma acyclicity, and cycle-freedom have two nice closure properties
each, which is optimal, by \thref{thm:limclosure}.
In other words these four notions are \emph{good acyclicity notions} in the sense
of the informal definition below.

But are there any other \emph{good acyclicity notions},
i.e. reasonable hypergraph acyclicity notions
that would also have two closure properties?
For example, can we find a reasonable property that is, like beta acyclicity,
closed under vertex and edge removal, but which is more general?
The answer is no: 
\begin{thmen}[optimality of the notions]
\label{thm:opt}
\begin{enumerate}
\item Gamma acyclicity implies any reasonable hypergraph acyclicity notion.
\item Beta acyclicity generalizes any reasonable hypergraph acyclicity notion that is closed under both edge removal and vertex removal.
\item Alpha acyclicity is the closure of gamma acyclicity for edge contraction; therefore it is the most restrictive reasonable hypergraph acyclicity notion that is closed under edge contraction.
\item Cycle-freedom generalizes any reasonable hypergraph acyclicity property that is closed under vertex removal.
\end{enumerate}
\begin{proofennoend}[1] Let $P$ be a reasonable property, and $\Hy$ a gamma acyclic hypergraph.
By \thref{charDM}, $\Hy$ is DM-reducible, i.e. we can obtain the empty hypergraph by applying DM operations.
Since $P$ is reasonable, $P$ is invariant w.r.t. any of the DM operations, by properties (1), (3), and (4) of reasonable notions.
Therefore, we have $P(\Hy)$ iff $P(\emptyset)$, which is true by property (2) of reasonable notions.
We have proved that gamma acyclicity implies $P$.
\end{proofennoend}
\begin{proofennoend}[2, sketch]
By \thref{lem:beta1b}, a hypergraph is beta acyclic if and only if every subset of it is cycle-free.
\end{proofennoend}
\begin{proofennoend}[3]
We already know by \thref{thm:contraction} that alpha acyclicity is closed under edge contraction.
We therefore only have to prove that any alpha acyclic hypergraph can be obtained by edge contraction of some
gamma acyclic hypergraph.

Let $\Hy$ be an alpha acyclic hypergraph.
Consider a join tree $(\mathcal{T},l)$ of $\Hy$.
For short, we say two edges $e$ and $f$ of $\Hy$
are neighbours in $\mathcal{T}$ if the vertices $a$ and $b$ of $\Hy$
labeled resp. $e$ and $f$ (i.e. $l(a)=e$ and $l(b)=f$) are neighbours in $\mathcal{T}$.
We define, for every $e\in\Hy$:
\mydisplaystyle{\begin{align*}
&f(e)=\left\{\left.x_i^{\{e,f\}}\,\right|\,x_i\in e\cap f\ \mathrm{and}\ f\ \mathrm{neighbour\ of}\ e\ \mathrm{in}\ \mathcal{T}\right\}
\end{align*}}
Let $\Hy'=\{f(e)\,|\,e\in\Hy\}$.
We prove $\Hy'$ is gamma acyclic.
One way to prove it consists in using Duris' characterization, stated \thref{char:gamma}:
notice $\Hy'$ has trivially a join tree with disjoint branches for any root, which establishes the result.
Another way to prove it consists in noticing that $\mathcal{D}(\Hy')$ is a usual tree.
As a consequence, it is gamma acyclic by \thref{rem:acygraph}, therefore, by \thref{thm:self-dual}, 
so is $\Hy'$.
We have proved that $\Hy'$ is gamma acyclic.

For each \smash{$x_i\in\Vertexp{\Hy}$} the vertices in the form \smash{$x_i^{\{a,b\}}$} are connected by join property, 
we proceed to edge contraction of all into a single vertex \smash{$x_i$}, and get the hypergraph $\Hy$.
This concludes the proof.
\end{proofennoend}
\begin{proofen}[4]
Easy from the definitions.
See appendix for details.
\end{proofen}
\end{thmen}
In other words, a \emph{good acyclicity notion} is either between
gamma and beta acyclicity, or between alpha acyclicity and cycle-freedom.
In particular, reasonable hypergraph acyclicity notions between beta and alpha acyclicity are, in the best case, closed under vertex removal.
If we focus on the properties closed under vertex removal --- the only properties that have a chance to have \emph{two} closure properties,
we get Figure~\ref{fig:closure}.
\begin{figure}~
\begin{center}
\mbox{\input{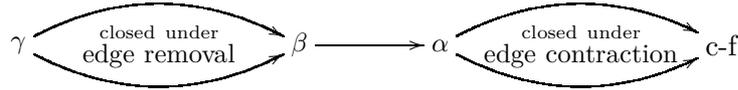}}
\end{center}
\caption{Classification of good acyclicity notions, i.e. reasonable acyclicity notions that are closed under vertex removal and either edge deletion or edge contraction. Here $\alpha$ (resp. $\beta$, $\gamma$, c-f) denote alpha acyclicity (resp.
beta acyclicity, gamma acyclicity, cycle freedom).\label{fig:closure}}
\vspace*{6pt}
\end{figure}
\begin{remen}
For any good acyclicity notion $P$ that is between gamma and beta acyclicity,
its dual, denoted $\mathrm{dual}(P)$, and defined as $\mathrm{dual}(P)(\Hy)\Leftrightarrow P(\mathrm{dual}(\Hy))$,
is also a good acyclicity notion. 
\end{remen}
As an example, Duris \cite{betagamma} introduced an acyclicity notion,
namely the fact of having a rooted join tree with disjoint branches,
that has found an application in \cite{CapelliDurandMengel}.
It is easy to prove this notion is a reasonable acyclicity notion,
and that it is closed under vertex and edge deletion.
Consequently, it is a good acyclicity notion, and, as stated by Duris,
it is between gamma and beta acyclicity.
Nevertheless, this notion does not have all the enjoyable properties of beta and gamma acyclicity:
it is not self-dual. But its dual is also a good acyclicity notion, by previous remark.
\mysection{Conclusion}
\noindent
We have proposed a consistent combined presentation of different characterizations of alpha, beta and gamma acyclicity.
The proposed characterizations cover the main use of these notion in some domains, such as finite model theory.
Moreover, we provided a ``proof framework'' that allows easy proofs of other characterizations.

Thanks to these different characterizations, we were able to prove that any 
``good acyclicity notion'' is either between gamma and beta acyclicity, or between alpha acyclicity
and cycle-freedom. 
Since alpha and beta acyclicity are on the two extremities of the ``gap'' between them (there is no good notion in between)
we conclude that these two notions are of particular interest.
\mysection{Acknowledgements}
\noindent
The author would like to thank \'Etienne Grandjean for careful proof reading,
Luc Segoufin for useful comments and suggestions, 
and Stefan Mengel for detailed corrections on the latter versions of this paper and 
for numerous discussions that directly motivated and inspired section 4.2,
in particular the main theorem of this section, \thref{thm:opt}.
%
%
%
%

\renewcommand{\bibsection}{\mysection{References}}
\bibliography{myrefs}		
\mysection{Appendix}
\begin{proofennoend}[\thref{thm:jt}]
We proceed by induction on the size of $\Hy$.
The empty tree is a join tree of the empty hypergraph.
Assume the equivalence holds for any hypergraph of size less than $n$.

Take $\Hy$ an alpha acyclic non-empty hypergraph of size $n$.
It is GYO-reducible, consider the first operation.
If it is an included edge removal, then the hypergraph obtained, called $\Hy'$, is also reducible hence acyclic,
therefore by induction it has a join tree $(\mathcal{T},L)$.
Consider the edge $e\in\Hy\smallsetminus\Hy'$.
It is included in some other edge $f\in\Hy'$.
Consider the vertex $x$ in $\mathcal{T}$ such that $L(x)=f$.
Take any fresh symbol $y\notin\Vertexp{\mathcal{T}}$.
Extend the definition of $L$ with $L':y\mapsto e$ and $L'(t)=L(t)$ for every $t\neq y$.
It is easy to check that \smash{$\big(\mathcal{T}\cup\{\{x,y\}\},L'\big)$}
satisfies the join property, and is therefore a join tree of $\Hy$.

In the other case, the first operation of the reduction is the removal of a singleton vertex $t$,
then $\Hy'=\Hy[\smallsetminus\{t\}]$ is GYO-reducible hence acyclic,
so by induction it has a join tree $(\mathcal{T},L)$.
Since $t$ is a singleton vertex in $\Hy$, there is only one edge $e_t$ that includes
it. The set $e_t\smallsetminus\{t\}$ is an edge of $\Hy'$,
we call $x$ the vertex of $\mathcal{T}$ such that $L(x)=e_t\smallsetminus\{t\}$.
Two cases: either $e_t\smallsetminus\{t\}$ belongs to $\Hy$ or not.
In the first case, take any fresh vertex $y\notin\Vertexp{\mathcal{T}}$,
define $L':y\mapsto e_t$ and $L'(z)=L(z)$ for any vertex of $\mathcal{T}$ other than $y$.
It is easy to check that \smash{$\big(\mathcal{T}\cup\{\{x,y\}\},L'\big)$}
satisfies the join property, and is therefore a join tree of $\Hy$.

In the other case, just define $L'(x)=e_t$ and $L'(z)=L(z)$ for any other vertex $z$,
then $(\mathcal{T},L')$ is a join tree of $\Hy$.

We have proved that an alpha acyclic hypergraph of size $n$ has a join tree, now we prove the converse.
Assume $\Hy$ has a join tree $(\mathcal{T},L)$, and consider a leaf $x$ of this join tree.
If $x$ is the only vertex of the join tree, then the result is obvious, assume this is not the case:
let $y$ be its neighbour.
We have two cases, either $L(x)\subset L(y)$ or not.

Assume $L(x)\!\subset\!L(y)$. Then $\Hy'=\Hy\smallsetminus\{L(x)\}$ has an obvious join tree,
therefore, by induction, $\Hy'$ is alpha acyclic so it is GYO-reducible, therefore so is $\Hy$,
therefore $\Hy$ is acyclic.

In the other case, take $t\in L(x)\smallsetminus L(y)$. By the join property,
$t$ is contained only in $L(x)$. If $L(x)\smallsetminus L(y)=\{x\}$, then 
\smash{$\big(\mathcal{T}[\smallsetminus\{x\}],L\big)$}
is a join tree of \smash{$\Hy\big[\smallsetminus\{t\}\big]$}; in the other case, defining $L':z\mapsto L(z)\smallsetminus\{t\}$
shows that \smash{$\Hy\big[\smallsetminus\{t\}\big]$} also has a join tree.
In both cases, \smash{$\Hy\big[\smallsetminus\{t\}\big]$} is alpha acyclic hence has an alpha elimination
order, and since $t$ is an alpha leaf of $\Hy$, $\Hy$ has an alpha elimination order
therefore it is acyclic.
\end{proofennoend}
\begin{proofennoend}[\thref{thm:opt} (2)]
Let $P$ be a reasonable hypergraph acyclicity notion, that is closed under edge removal and vertex removal.
Assume $P$ is not a particular case of beta acyclicity; there must therefore be some $\Hy$ such
that $P(\Hy)$ but $\Hy$ is not beta acyclic.
We can get a usual graph cycle from $\Hy$ be edge and/or vertex removal, therefore some usual graph
cycle satisfies $P$, which is therefore not a reasonable acyclicity notion.
\end{proofennoend}
\begin{proofennoend}[\thref{thm:opt} (4)]
Let $P$ be a reasonable acyclicity notion that is closed under vertex removal,
and let $\Hy$ be a hypergraph.
Suppose $\Hy$ is not cycle-free. Then we can find $S$ such that $\Hy[S]$ is a usual graph cycle.
By property (2), $P(\Hy[S])$ does not hold hence, since $P$ is closed under vertex removal, $P(\Hy)$ does not hold.
We have proved that a hypergraph satisfies $P$ only if it is cycle-free.
\end{proofennoend}

\NumT
\end{document}